\input amstex.tex

\input amsppt.sty

\TagsAsMath

\NoRunningHeads

\magnification=1200

\hsize=5.0in\vsize=7.0in

\hoffset=0.2in\voffset=0cm

\nonstopmode

\document

\document

\topmatter

\title{Dispersion   for Schr\"odinger equation with
periodic potential in 1D}
\endtitle

\author
Scipio Cuccagna
\endauthor

\address
DISMI University of Modena and Reggio Emilia, Via Amendola 2,
Padiglione Morselli, Reggio Emilia 42100 Italy
\endaddress

\email cuccagna.scipio\@unimore.it\endemail

\abstract We extend a result on dispersion for solutions of the
linear Schr\"odinger equation, proved by Firsova for operators with
finitely many energy bands only, to the case of smooth potentials in
1D with infinitely many bands. The proof consists in an application
of the method of stationary phase. Estimates for the phases,
essentially the band functions, follow from work by Korotyaev. Most
of the paper is devoted to bounds for   the Bloch functions. For
these bounds we need   a detailed analysis of the quasimomentum
function and   the uniformization of the inverse of the
quasimomentum function.

\endabstract

\thanks Research fully supported by a special grant from the Italian
Ministry of Education, University and Research.
\endthanks

\endtopmatter

\head \S 1 Introduction \endhead

We consider  an operator $H_0= -\frac {d^2}{dx^2} +P(x)$ with $P(x)$
a smooth periodic potential   of period 1.  We will prove:
\proclaim{Main Theorem  } There is a $C>0$ such that for any $p\ge
2$ and   $t>0$ we have $ \| e^{itH_0} \colon L^{\frac{p}{p-1} }(\Bbb
R)  \to  L^{p }(\Bbb R) \| \le C \max \left  \{ t^{-\frac 12},  t
^{-\frac 13} \right \} ^{(1-\frac 2p)} . $
\endproclaim
 We recall that the spectrum $\Sigma
(H_0)$   is a union of closed intervals (bands) and that the theorem
in the case of finitely many bands is in Firsova \cite{F1}, for
alternative proofs in the two bands case see \cite{Cai,Cu}.
 Here we will consider the case when $\Sigma (H_0)$ is a
union of infinitely many bands. Just to simplify the notation, we
will formulate the proof only in the generic case when all the gaps
are nonempty. This generic case contains all the essential
difficulties, and a proof for cases when some gaps are empty goes
through similarly, only with    more complicated notation. We are
motivated
 by nonlinear problems, see \cite{Cu}.
Indeed dispersion for linear operators is central to nonlinear
equations, see for instance \cite{Str}. For $P\equiv 0$ we have $
e^{-it {d^2}/{dx^2}} (x,y)= (4\pi it)^{-\frac{1}{2}} e^ {-i
{(x-y)^2}/{4t} }$. For no nonconstant $P(x)$ a similar explicit
formula seems to be available.
 We outline the proof, with terminology and  formulas  introduced
rigourously later. The  proof of the Main Theorem reduces to a
pointwise
  bound on the kernel $e^{itH_0}(x,y).$ The
kernel is expressed   by means of the distorted Fourier transform,
    written using Bloch functions and is  a  sum of
oscillatory integrals, one term  per energy band. In each integral
the phase involves the band functions, which express energy $E$ with
respect to quasimomentum $k$. $k$ varies in $\Bbb C$ cut along some
slits which correspond to the gaps. We derive various estimates for
the derivatives of the band function $E(k)$ from the representation
of the quasimomentum function $k=k(w)$,   $w=\sqrt{E}$, exploiting
formulas in Korotyaev \cite{K1}. We have essentially $E(k)=k^2$,
except for very narrow regions near the edges of the spectral bands
where, in particular,
   the third derivative of $  E(k)$
 is very large. In this region very close the edges our bounds on the Bloch functions are large, but in the stationary phase
  formula are more than offset by the
  very small upper bound  for the inverse  of the third derivative of $
  E(k)$. Away from the edges, the contribution  to $e^{itH_0}(x,y) $
  is essentially the same of the constant coefficients case and
  there are good bounds for the Bloch functions. There is an
  intermediate region, close  but not very close, in a relative sense,   to the edges, where our bounds for the Bloch
  functions are large and where the third derivative of $  E(k)$
  is not large. Yet,  combining the narrowness of this region with
  good enough bounds for the Bloch
  functions,   we control the corresponding contribution to $e^{itH_0}(x,y) $.
Since estimates for $E(k)$ follow  by fairly direct  elaboration on
material by Korotyaev
   \cite{K1}, our most serious effort is for
   bounds on the Bloch functions and their first
  derivative in $k$.  It is probably because of these  that \cite{F1}  considers
   finitely many bands only. The Bloch functions
  can be expressed as a linear combination of  a nice fundamental
  set of solutions of the equation $H_0u=Eu$. The coefficient in
  this linear combination is the Weyl-Titchmarsh function, which
  near the   edges of the bands of higher energy    is hard to bound, because it is the ratio
  of two very small quantities. In the finite
  bands case, Firsova
\cite{F1} uses the fact that the Bloch functions are analytic on the
uniformization of the function $w(k)$. Since there are only finitely
many bands, this gives a uniform bound. Bounds for the derivatives
follow from the Cauchy integral formula. In the case of infinitely
many bands,  near each edge the Bloch functions are
  bounded, but there is no obvious uniform bound over an
infinite number of bands. Furthermore, the distance of the edges
from the boundary of the domain of analyticity goes to 0 as we take
higher energies. Hence,  for large energy, the Cauchy integral
formula gives bad estimates on the $k$ derivative of the Bloch
functions, near the edges of the band. This problem for the
derivative is the main source of trouble in the paper. Taking
Fourier expansion, if one differentiates in $k$ one gets  a small
divisors problem, with small divisors for   modes $n=0$ and for $\pm
n\pi +k\approx 0$. Is seems convenient to use the Cauchy integral
formula with bounds for complex $k$. The Bloch functions are of the
form $e^{\pm ikx}m_{\pm }(x,k)$. Thanks to estimates on a nice pair
of fundamental solutions of $H_0u=Eu$ and relations between the
various terms in the Weyl-Titchmarsh function and a certain
normalization term denoted by $N(k)$, we see that for large energies
in  the Fourier series expansion $m_{\pm }(x,k)=\sum \widehat{
m}_{\pm }(n,k)e^{2\pi i nx}$ most terms are small compared to
$\widehat{ m}_{\pm}(0,k)$ and $\widehat{ m}_{\pm}(\pm n,k)$ with $|
\pi n +k|\ll 1$ and  $\| m_{\pm }(\cdot ,k)\| _{L^\infty}\approx \|
m_{\pm }(\cdot ,k)\| _{L^2}$. For $k\in \Bbb R$ the Bloch functions
are normalized, $\| m_{\pm }(\cdot ,k)\| _{L^2}=1$. However for
$k\not \in \Bbb R$ and near the boundary of the domain of
analyticity, that is near the slits,  it is problematic to bound $\|
m_{\pm}(\cdot ,k)\| _{L^2}.$
 In correspondence to the interior of the band   $|\widehat{ m}_{\pm }(0,k)|\gg |\widehat{ m}_{\pm }(\mp n,k)|$.  Using
the normalization of the Bloch functions we get
$$1\approx \widehat{ m}_{+}(0,k)\widehat{ m}_{-}(0,k)+\widehat{
m}_{+}(-n,k)\widehat{ m}_{-}(n,k),\tag 1.1$$  so
 $\| m_{+ }(\cdot ,k)\| _{L^2} \| m_{-}(\cdot ,k)\| _{L^2} \lesssim 1$.
Near the  slits however $|\widehat{ m}_{\pm }(0,k)|\approx
|\widehat{ m}_{\pm }(\mp n,k)|$ and     in the right in (1.1) we
could have a cancelation. However from some explicit formula and
thanks to (9.1) we get $\widehat{ m}_{+}(0,k)\widehat{
m}_{-}(0,k)-\widehat{ m}_{+}(-n,k)\widehat{ m}_{-}(n,k)\approx \dot
E/(2k)$. By the Schwartz Christoffel formula we get a good bound for
this quantity for $k$ close enough to the real axis. Hence  for $k$
close enough to the real axis in (1.1) there is no cancelation and
we get a uniform bound on $\| m_{+ }(\cdot ,k)\| _{L^2} \|
m_{-}(\cdot ,k)\| _{L^2} $. The argument shows however that near the
extremes of the slits the two terms in (1.1) are in fact unbounded.
Away from the edges of the bands we improve the estimates
significantly thanks to more information on $E(k)$.

   One can relax  significantly the regularity
requirements on $P(x)$ maintaining the   proof.  The proof goes from
\S 4 to \S 10. In particular, the estimates on the Bloch functions
are in \S 9 and \S 10.

Here the spectrum is denoted by $\Sigma (H_0)$, with  the usual
notation $\sigma (H_0)$  reserved for something else. In a statement
or in a proof, notation $a_n\lesssim b _n$ means that there is a
fixed constant $C>0$ independent of $n$, with $a_n\le C b _n$. If we
write $a_n\ll b _n$ we mean $a_n\le C b _n$ for a very small fixed
constant $C>0$. If we write $a_n\approx b _n$ we mean $(1/C)b_n \le
a_n\le C b _n$ for a fixed constant $C>0$ independent of $n$. If we
write $a_n=O(b_n)$ (resp. $a_n=o(b_n)$) we mean $a_n\lesssim b _n$
(resp. $a_n\ll b_n$). For $p\in [1,\infty ]$, by $\| f \|_p$ we mean
the usual $L^p$ norm of $f(x)$, where the $x$ varies in a set
indicated in the context.

\head \S 2 Band function, Bloch function, quasimomentum and
uniformization
\endhead

The spectrum is of the form $\Sigma (H_0)=\cup _{n=0}^\infty \Sigma
_n$, with each two compact intervals $\Sigma _n$  and $\Sigma
_{n+1}$ separated by an open gap $G_n$. We assume  $\inf \Sigma
_0=0.$ We set now $ \sigma  = \cup _{n=-\infty}^\infty \sigma _n$,
with each two compact intervals $\sigma _n$  and $\sigma _{n+1}$
separated by  an open gap $g_n$, with   $g_0$ empty, and with
$\sigma _{-n}=-\sigma _n$  and $\sigma ^2_n=\Sigma _n$. We will
assume each $g_n$ non empty for $n\neq 0$. We recall now the
following standard result, see \cite{Ea} ch. 4: \proclaim{Theorem
2.1} Let $P(x)$ be smooth. Set $\sigma _n=[a_n^+, a_{n+1}^-]$ and
$g_n=]a_n^-, a_{n }^+[$. Then there exist a strictly increasing
  sequence $\{ \ell _n \in \Bbb Z \} _{n\in \Bbb Z}$ and a
fixed constant $C $ such that
$$|a_n^- -\ell _n\pi |+ |a_n^+ -\ell _n\pi |\le C  \langle \ell _n \rangle ^{-1}.
$$
For any $N$ there exists a fixed constant $C_N $   such that the
length $|g_n|$ of the gap $g_n$ is $|g_n|\le C _N \langle \ell _n
\rangle ^{-N}. $
\endproclaim
To simplify notation we will assume in the rest of the paper that
$\ell _n\equiv n$, the case when all spectral gaps, from a certain
one on,  are not empty, which is generic, but  essentially the same
proof goes through in general. For any $w  \in {\Bbb C_+}$ (the open
upper half plane)  we consider the fundamental
   solutions $\theta (x,w)$ and
$\varphi (x,w)$ of $H_0u=w^2u$ which satisfy  the initial conditions
$$\varphi (0,w)=   \theta ^\prime (0,w) =0\, , \quad               \varphi ^\prime (0,w) =\theta
(0,w)=1.\tag 2.1$$ The Floquet  determinant $D(w)$ is defined by
$2D(w)={\varphi ^\prime (w) +\theta (w)}$ where $\varphi ^\prime (w)
=\varphi ^\prime (1,w)$ and   $\theta  (w) =\theta (1,w)$. For any
$w  \in  {\Bbb C_+}$   there is a unique   $  k\in {\Bbb C_+}$,
called quasimomentum, and a unique choice of constants $m^\pm (w)$
such that the functions
$$\tilde \phi _\pm (x,w)=\theta (x,w)+m^\pm (w)
\varphi (x,w)\tag 2.2$$ are of the form $\tilde \phi _\pm
(x,w)=e^{\pm ikx} \xi _\pm (x,w)$ with  $\xi _\pm (x,w)$ periodic of
period $1$ in $x$. We have
$$m^\pm (w)=\frac{\varphi ^\prime (w)  - \theta (w) }{2\varphi (w)}
\pm i \frac{\sin k}{\varphi  (w)}.\tag 2.3$$ We have the relation
$D(w)=\cos k$. The correspondence between $w$ and the corresponding
quasimomentum $k$ is a conformal mapping between ${\Bbb C_+}$  and a
"comb" $K$, that is a set $K={\Bbb C_+}\backslash \cup _{n\neq 0}  [
n\pi , n\pi +ih _n ]$ where the $[ n\pi , n\pi +ih _n ]$ are
vertical slits with $h_n\ge 0$. In particular, $|g_n|\le 2 h_n \le
(1+C n^{-2}) |g_n|$ for a fixed $C$, see Theorem 1.2 \cite{KK}. Now
we will use that all gaps are nonempty, but the following standard
discussion extends easily.
 The map $k(w)$ is called quasimomentum map and
extends into a continuous map in $\overline{\Bbb C_+}$ with
$k(\sigma _n)=[ n\pi , (n +1)\pi ]$, with $k(w)$ a one to one and
onto map between $\sigma _n$ and $[ n\pi , (n +1)\pi ]$, and with
 $k(g _n)=] n\pi ,
n\pi +ih_n ]$. We have $k(-\bar w) =- \overline{k(w)}$. By the
Schwartz reflection principle,   $k(w)$ extends into a conformal map
from $\Bbb C \backslash \cup _{n\neq 0}\overline{g _n}$ into $\Cal
K=\Bbb C\backslash \cup _n \gamma _n $  with $\gamma _n= [n\pi
-ih_n, n\pi +ih_n]$. So we have $k( \bar w) = \overline{k(w)}$  and
$k(w)=-k(-w)$. Hence also $ w( \bar k)=  \overline{w(k)}$ and
$w(k)=-w(-k)$. Then for $|t|<h_n$ we have $w(n\pi +it \pm  0)=
-\overline{w(-n\pi +it \mp
 0)}= -w(-n\pi +it \mp  0)$. This means that the band function
$E(k)=w^2(k)$ extends in an analytic map with values in $\Bbb C$ and
with domain the Riemann surface $\Cal R$ obtained  identifying $n\pi
+it \pm  0$ and $-n\pi +it \mp  0$ for each $n$ and for each
$|t|<h_n$.

For $w\in \Bbb C^+$ we have introduced   $\tilde \phi _\pm
(x,w)=e^{\pm ik(w)x} \xi _\pm (x,w)$ with $\xi _\pm (x,w)$ periodic.
These functions extend by continuity to  $w\in \Bbb R$. For $w\in
\Bbb C^+$,
  we have  by the properties of $k(w)$ and by the definition of $\overline{\tilde \phi
_\pm (x,w)} $,  $$\overline{\tilde \phi _\pm (x,w)} =e^{\mp
i\overline{k(w)}x}\overline{\xi _\pm (x,w)}= e^{\pm i
k(-\overline{w}) x}\overline{\xi _\pm (x,w)}= \tilde \phi _\pm
(x,-\overline{w})  \tag 2.4 $$ and, for $w\in \sigma$, we have, for
the same reasons,
$$\align & \overline{\tilde \phi _\pm (x,w)} =e^{\mp i {k(\overline{w})}x}\overline{\xi
_\pm (x,w)}= e^{\mp i {k(\overline{w})}x} {\xi _\mp
(x,\overline{w})}= \tilde \phi _\mp (x, \overline{w})  ,  \tag
2.5\\& \tilde \phi _\pm (x,w) =e^{\pm i {k(w)}x} \xi _\pm
(x,w)=e^{\mp i {k(-w)}x} \xi _\pm (x,w) =\tilde \phi _\mp (x,-w) .
\tag 2.6 \endalign$$ By (2.6) the function $\tilde \phi _\pm (x,w)$
can be extended across $\sigma $ into analytic functions in $w\in
\Bbb C \backslash \cup _{n\neq 0}\overline{g _n}$ setting $\tilde
\phi _\pm (x,w)= \tilde \phi _\mp (x,-w).$ It is elementary to see
that (2.5) and (2.6) are now true for any $k\in \Cal K=\Bbb C
\backslash \cup _{n\neq 0} \gamma _n$. Set now
$$N^2(w)= \int_0^{1} \tilde \phi _+ (x,w)\tilde \phi _- (x,w) dx.$$
By (2.5) we have $N^2(w)= \int_0^{1}   \big |\tilde  \phi _\pm (x,w)
\big |^2 dx>0$ for $w\in \sigma $ (so that we define $N(w)>0$ for
$w\in \sigma $).    $N^2(w)$ is well defined and analytic in $\Bbb C
- \cup _{\neq 0}g _n$. We  have
$N^2(w)=N^2(-w)=\overline{N^2(\overline{w})}$, the first equality by
(2.6) and the second by (2.5). From formula (1.4) \cite{F2} we have
$D' (w)=-4w \varphi (w)N^2(w)$, for a sketch of proof see \S 3
\cite{Cu}. From $D(w)=\cos k$ we get $D' (w)= - \frac{dk}{dw}\sin
k$. Since $k(w)$ is a conformal map,  $\frac{dk}{dw} \neq 0$ for
$\Im w
>0$ and so for $k\in K$. Hence $N^2(w)\neq 0$ for any $w\in \Bbb
C \backslash  \cup _{n\neq 0}\overline{g _n}$.   We set now
$$ e^{ik (x-y)}m _+ ^0(x,w)m _- ^0 (y,w)=\frac{\tilde \phi _+(x,w)\tilde \phi _-(y,w)}{N^2 (w)} .\tag 2.7$$ We
express $w=w(k)$ for $k\in K$ and with an abuse of notation we write
   $m^0 _\pm (x,k)$ for
$m^0 _\pm (x,w(k))$. Then the product $m _+ ^0(x,k)m _- ^0 (y,k)$
extends analytically for $k\in \Cal K$ and  to $\Cal R$. Fix a
square root $N$ of $N^2$ with $N(k)>0$ for $k\in \Bbb R$, ignoring
problems of monodromy (here we need the product $m _+ ^0 m _- ^0$ to
be analytic). We set $m _\pm ^0(x,k)=\xi _\pm ^0(x,k)/N(k)$.

\bigskip

\head \S 3 Fourier transform  \endhead

From Theorem XIII.98 \cite{RS} it is possible to conclude:

\proclaim{Lemma 3.1} Let $N(k)=\sqrt{N^2(k)}>0$ for $k\in \Bbb R$.
Set $\phi  _\pm  (y,k)=\widetilde{\phi}  _\pm  (y,k)/N(k)$ and $\hat
f (k)=\int _\Bbb R {\phi  _+ (y,k)} f(y)dy$. Then:

$$\align & \int  _\Bbb R| f(y)|^2dy=  \int  _{\Bbb R}|\hat f (k)|^2dk \tag a
\\& f(x)=   \int _{\Bbb R}  {\phi _-
(x,k)} \hat f(k) dk \tag b
\\& \widehat{H_0f  }(k)= E (k)\hat f(k) . \tag c \endalign
$$
\endproclaim

Lemma 3.1 implies:

 \proclaim{Lemma 3.2} We have
$  e^{itH_0} (x,y)=      K (t,x,y)  =\sum _{n\in \Bbb Z} K^n(t,x,y)$
with

 $$ K^n
(t,x,y)=
  \int   _{n\pi }^{(n+1) \pi}
e^{i(t E  (k) -(x-y)k) }   {m _- ^0(x,k)}{ m  _+ ^0(y,k)} dk . \tag
3.1
$$
\endproclaim

Then the Main Theorem   follows from:

\proclaim{Theorem  3.3} There is $C$ fixed such that $
 | K(t,x,y)|\le C \max \{     t   ^{-\frac 13}, t  ^{-\frac 12} \}  .
 $
\endproclaim

 \noindent  Theorem  3.3 follows by  the method of stationary phase. Notice that the phase in
$K^n(t,x,y)$ satisfies the following   result by Korotyaev
\cite{K2}:

\proclaim{Theorem   3.4} Consider    $E(k)$   for $k\in [ \pi n, \pi
(n+1)]$. Then $ E   ' (k)=0$ for $k=n\pi , (n+1)\pi  $ and $ E '
(k)>0 $ in $]\pi n, \pi (n+1) [$ for $n\ge 0$ ($E(k)$ is even). In
$[ \pi n, \pi (n+1)] $ the equation $E   '' (k)=0$ admits exactly
one solution $k_n$. We have $k_n\in ] \pi n, \pi (n+1) [$ and $E
^{\prime \prime \prime } (k_n)\neq 0$.
\endproclaim

Naively the proof of Theorem 3.3 would go as follows.
 Theorem 3.4 and estimates on   $m _\pm ^0(x,k)$
for  $k\in [ \pi n, \pi (n+1)]$ lead to an estimate $|K^n(t,x,y)|\le
D_n  \langle  t \rangle ^{-\frac 13}$ thanks to the method of
stationary phase, which we quote   from p. 334 \cite{Ste}:

\proclaim{Lemma  3.5} Suppose $\phi (x)$ is real valued and smooth
in $[a,b]$ with $|\phi ^{(m)}(x) |\ge c_m>0$ in $]a,b[$ for $m\ge 1$
. For $m=1$ assume furthermore that $\phi ^\prime (x)$ is monotonic
in $]a,b[$. Then   we have for $C_m=  5 \cdot 2^{m-1}  -2 $:
$$ \big |\int _a^b   e^{i\mu \phi (x)}  \psi (x)  dx             \big
| \le C_m  (c _m\mu )^{-\frac 1m} \left [ \min \{ |\psi (a) |, |\psi
(b) | \}+\int _a^b |\psi ' (x) | dx \right ] . $$
\endproclaim

Since we need     to add up over   all the $K^n$   we have to
control the constants $D_n$. In the next section we state   a list
of estimates on the band function $E(k)$ and on the Bloch functions
$\phi _\pm (x,k)$ which are sufficient to the purpose of the present
paper and  then prove Theorem 3.3. In the subsequent sections we
prove the estimates on $E(k)$ and $\phi _\pm (x,k)$.

\bigskip

\head \S 4 Proof of Theorem 3.3 \endhead

We state estimates for the first three derivatives of $E(k)$. These
are proved later using material in Korotyaev \cite{K1}. Notice that
$E(k)$ is even in $k$, so we consider only $k\ge 0$. We start with
the first derivative $\dot E (k)$:

\proclaim{Lemma    4.1} For all $k\ge 0$ we have $\dot E(k)  \ge 0$.
  $\dot E(k) =0$ implies $k=n\pi $ for some $n\in \Bbb Z$. There are  fixed constants
$C>0$ and $c>0$ such that for any $n$ we have
$$\aligned &   \frac{k}{C}   \frac {\sqrt{
 w-a^+_n    }    }{\sqrt{|g_n|}  } \le \dot E(k) \le C k   \frac {\sqrt{
 w-a^+_n    }    }{\sqrt{|g_n|}  }  \quad \text{for} \quad a^+_n \le
 w
 \le a^+_n + c|g_n| \\& |\dot E(k) - 2k| \le \frac{C}{\langle k \rangle
  } \quad \text{for} \quad   a^+_n + c|g_n| \le w   \le a^-_{n+1}- c|g_{n+1}|
  \\& \frac{k}{C}\frac {\sqrt{
  a^-_{n+1}  -w  }    }{\sqrt{|g_{n+1}|}  } \le
  \dot E(k) \le  k C \frac {\sqrt{
  a^-_{n+1}  -w  }    }{\sqrt{|g_{n+1}|}  }  \quad \text{for} \quad  a^-_{n+1} - c|g_{n+1}|\le
 w
 \le  a^-_{n+1} .\endaligned $$

\endproclaim

 Now we consider the second
derivative $\ddot E (k)$:

 \proclaim{Lemma 4.2} There are fixed constants
$C>0,$ $C_1> C_2>0$ and $c>0$ such that for  any $n$  and any $k\in
[n\pi , (n+1)\pi ] $, that is for any $w\in [a^+_n , a^-_{n+1}]  $,
we have:

$$\aligned    &     w
 \le a^+_n  + c|g_n|   \Rightarrow              \big |\ddot{E}-\frac{n}{|g_{n  }|} \big
 | \le  C \\& w
 \ge a^-_{n+1}  - c|g_{n+1}|  \Rightarrow              \big |\ddot{E}+\frac{{n+1}}{|g_{{n+1}  }|} \big
 | \le  C\\& a^+_n +c|g_n| \le w\le     a^-_{n+1}-C _1|n+1|^{\frac{1}{3}}|g_{n+1
}|^{\frac{2}{3}}   \Rightarrow             \ddot{E} \approx \frac 12
+\frac{n|g_n|^2}{ |w-a_n^+|^{3}} .
 \endaligned $$
\endproclaim

Finally we have (see Lemma 7.5):

\proclaim{Lemma 4.3} There are   $c>0$ and $c_1>0$ such that
$\forall$ $\alpha \in (1/2, 1)$  and  for any $w\in \sigma _{n}\cup
\sigma _{n-1}$ with $c|g_{n }|\le |w - a^\pm _{n }| \le |g_{n }|^{
\alpha}$ we have $ |\dddot{E} |\ge c_1 |n| \, |g_n|^{4(\frac{1}{2}-
\alpha )}.$

\endproclaim

Next we
 need estimates for the Bloch functions. First of all we have, see in  \S 9:

\proclaim{Lemma  4.4} There  are fixed constants $C>0$, $C_3>0$,
$\delta >0$, $\Gamma >0$ and $c>0$ such that for all $x$, all $n$ we
have :

{\item{ (1)}}   $  \forall \,   w \in [a^+_{n }
 +C_3n^5 |g_n| ,a^-_{n+1 }-C_3(n+1)^5 |g _{n+1}| ]$ we have
$$\aligned & \big |  m ^0_+(x,k)  m ^0_-(y,k) -1 \big | \le    C \langle
k\rangle ^{-1};
\endaligned$$

{\item{ (2)}}  for all $k=p+iq$ with $   k\in \Cal K$, $\pi n < p<
\pi (n+1)$, $k$ in $\{ k: |q|<\delta |g_n|\} \cup \{ k: |\pi n -
p|>\Gamma  |g_n|\}       \cup      \{ k: |\pi (n+1) - p|>\Gamma
|g_{n+1}|\}   $ and $|q|\le 1$ we have
 $$\big |  m ^0_+(x,k)  m ^0_-(y,k)   \big | \le C  .
 $$

\endproclaim

The proof of (1) Lemma 4.4 is elementary and that of (2)  if $k\in
\Bbb R$, that is the case near the edges not covered by (1), is
relatively easy. We will state later more estimates for $k\not \in
\Bbb R$, needed to bound   $\dot m ^0_\pm (x,k)= {\partial _k} m
^0_\pm (x,k)$.  We have:

  \proclaim{Lemma  4.5} There  are fixed
constants $C>0$ and $C_4>0$, with $C_4<C_2$, $C_2$ the constant in
Lemma 4.3,  such that for all $x$, all $n$ and for $v=0$ we have :

for  $a^{+}_{n}+  |g_n|^{\frac{1}{4}}\le u \le
\frac{a^{+}_{n}+a^{-}_{n+1}}{2} $,   there is a $C$ such that for
the corresponding $k=p+i0$ we have $  | \partial _k (m_- ^0(x ,k)
m_+ ^0(y ,k))\big | \le \frac{C}{k |k-\pi n|}
 ;$

  if  $
 \frac{a^{+}_{n}+a^{-}_{n+1}}{2} \le u \le a^{-}_{n+1}-
 |g _{n+1}|^{\frac{3}{5}}$ then $  | \partial _k
(m_- ^0(x ,k)  m_+ ^0(y ,k)) \big | \le \frac{C}{k |k-\pi (n+1)|}
 ;$

for all $    w $ in the remaining part of $  [a^+_{n },  a^-_{n+1
}]$ we have for
 a fixed $ C$ and with $m=n$ (resp. $m=n+1$)  near $a^+_{n }$ (resp. $a^-_{n+1 }$)
 $$\big | \partial _k
(m_- ^0(x ,k)  m_+ ^0(y ,k))    \big | \le    C  \left (    |k-\pi
m|
 +|g_m| \right ) ^{-1}.
 $$

\endproclaim

\bigskip
We assume now the above lemmas and go ahead with Theorem 3.3.

\bigskip

{\it Proof of Theorem 3.3.} Recall  $ K^n(t,x,y)$ given by (3.1). By
the discussion in \S 3, the only interesting case is when $|n|\gg
1$. Here we will sum over $n\gg 1$, the proof of the $n\ll - 1$
being similar. We  fix  a constant    $c\gg 1$. For $n\gg 1$ we
partition
$$ \aligned &\left [a_{n}^+ , a^-_{n+1} \right ]= \left [ a_{n}^+ , a_{n}^+ + c|g_n|\right ]
\cup
 \left [ a_{n}^+ + c|g_n|,    a_{n}^+ +
  |g_n|^{ \frac{1}{4}} \right ]  \cup \\&   \left  [ a_{n}^+ +
|g_n|^{ \frac{1}{4}},  a^-_{n+1}-  |g _{n+1}|^{ \frac{3}{5}} \right
]   \cup   \left [ a^-_{n+1}-  |g _{n+1}|^{ \frac{3}{5}} ,
a^-_{n+1}-    c |g _{n+1}|\right ] \\& \cup   \left [ a^-_{n+1}- c
|g _{n+1}|,a^-_{n+1}\right ]    . \endaligned$$  Here $ c\gg 1$ so
that for $w\in [a_{n}^++c|g_n|, a^-_{n+1}-c|g_ {n+1}|]$   we have $
\dot w \approx 1$.   For $w\in [a_{n}^++c |g_n|^{ \frac{1}{4}} ,
a^-_{n+1}- |g _{n+1}|^{ \frac{3}{5}}|]$ we have $ \ddot E\approx 1$.
We introduce a smooth, even, compactly supported cutoff $\chi
_0(t)\in [0,1]$ with $\chi _0\equiv 1$ near 0 and $\chi _0\equiv 0$
for $t\ge 2/3$. Set $\chi _1 =1-\chi _0 $. We split each $K^n=\sum
_1^5K^n_\ell $ partitioning the identity in $\sigma _n=[a_{n}^+,
a^-_{n+1}]$
$$\aligned &  1_{\sigma _n}(w)= \chi _0 (\frac{w - a_{n}^+ }{c|g_n|}  ) +  \chi _1 (\frac{w - a_{n}^+ }{c|g_n|}  )
  \chi _0(\frac{w - a_{n}^+ }{
 |g_n|^{ \frac{1}{4}}} ) + \\&   +  \chi _1 (\frac{w
- a_{n}^+ }{  |g_n|^{ \frac{1}{4}}}  )   \chi _1 (\frac{ a_{n+1}^--w
}{ |g _{n+1}|^{ \frac{3}{5}}}  ) + \chi _0(\frac{ a_{n+1}^--w }{ |g
_{n+1}|^{ \frac{3}{5}}} )  \chi _1  (\frac{ a_{n+1}^--w }{c |g
_{n+1}|} ) +\chi _0(\frac{ a_{n+1}^--w }{ c|g _{n+1}|} ) .
\endaligned
$$
 We bound one by one the
$K^n_\ell $.

\proclaim{Claim $\ell =1$ and $\ell =5$} For any $\epsilon >0$ there
is a fixed $C_\epsilon $ such that $ |K^n_1|  \le C_\epsilon
t^{-\frac{1}{2}}| g_{n}|^{\frac{1}{2}- \epsilon}$ and $ |K^n_5|  \le
C_\epsilon t^{-\frac{1}{2}}| g_{n+1}|^{\frac{1}{2}- \epsilon}$.
\endproclaim
{\it Proof of $\ell =1$.} By Lemma 3.5, by $\ddot E \ge
\frac{|n|}{2| g_{n}|}$, Lemma 4.1,  and by Lemmas 4.4 and 4.5,

$$|K^n_1(t,x,y)| \le \frac{C\sqrt{| g_{n}|}}{\sqrt{|n| t}} \int
_{a^{+}_{n}} ^{a^{+}_{n}+c|g_n|} \frac{ \frac{dk}{dw} dw}{ |k-\pi n|
 +|g_n| }  .
$$
By Lemma 4.1, $\frac{dk}{dw}\approx \frac{\sqrt{|
g_{n}|}}{\sqrt{w-a^{+}_{n}}}$ and $|k-\pi n| \approx \sqrt{|g_n|}
\sqrt{w-a^+_n} $. Hence

$$|K^n_1(t,x,y)| \le \frac{C }{\sqrt{ |n|  t}} \int
_{a^{+}_{n}} ^{a^{+}_{n}+c|g_n|}  \frac{ dw}{\sqrt{w-a^{+}_{n}}} \le
\frac{C_ 1 \sqrt{|g_n|  } }{\sqrt{ |n| t}}.
$$
With a similar argument we get the estimate for $K_5$.

\bigskip

\proclaim{Claim $\ell =2$  } There is $C>0 $ such that $ |K^n_2| \le
C \min\{  t^{-\frac{1}{2}}\log (1/|g_n|),  | g_{n}|^{\frac{1}{4}}
\}$.\endproclaim {\it Proof.} For $k\in [a^{+}_{n}+ c|g_n|, a^+_{n
}+   |g_{n }|^{\frac{1}{4}}]$ we have $|\ddot E| \gtrsim 1$ by Lemma
4.2 and $\frac{dk}{dw}\approx 1$ and $w-a^+_n \approx k-\pi n$ by
Lemma 4.1. Hence $ |K^n_2(t,x,y)| \le $

$$\le \frac{C _1}{\sqrt{  t}} \int _{a^{+}_{n}+c|g_n|} ^{a^{+}_{n}+
|g_{n }|^{\frac{1}{4}}} \frac{ \frac{dk}{dw} dw}{ |k-\pi n|
 +|g_n| }  \le C    t^{-\frac{1}{2}}     \int _ {a^+_{n }+  c|g_{n }| }
^{ a^+_{n }+   |g_{n }|^{\frac{1}{4}}} \frac{ dw}{ w-a^+_{n }
  }  \le A t^{-\frac{1}{2}}
 \log \frac 1{|g_n|}  .
$$
On the other hand, taking absolute value in the integral defining
$K^n_2$, thanks to Lemma 4.4 we get $|K^n_2|\le C
|g_n|^{\frac{1}{4}}.$

 \bigskip \proclaim{Claim  $\ell =4$} There is $C  $ such that   $ |K^n_4|
\le C  t^{-\frac{1}{3}}  | g_{n+1}|^{  \frac 1 {31}}$.
\endproclaim
{\it Proof.} $K^n_4$ is   an integral on an interval where $|\dddot
E| \gtrsim |n+1|  |g_{n+1 }|^{-\frac{1}{10}}$ by Lemma 4.3.
$\frac{dk}{dw}\approx 1$ and $ a^-_{n+1}-w \approx \pi (n+1)-k$ by
Lemma 4.1. So, by Lemmas 3.5, 4.4 and 4.5, we get a contribution
bounded by
$$   t^{-\frac{1}{3}}  |n+1|^{\frac{1}{3}} |g_{n +1}|^{ \frac{1}{30}}  \int ^ {a^-_{n +1}- c|g_{n+1 }| }
_{ a^-_{n +1}- |g_{n+1 }|^{\frac{3}{5}}} \frac{ dw}{ a^-_{n+1}-w
  }  \lesssim   t^{-\frac{1}{3}}  |n+1|^{\frac{1 }{3} } |g_{n+1 }|^{ \frac{1}{30}}
 \log \frac 1{|g_{n+1 }|}  .
$$
\bigskip

We get $ \sum _n \sum _{\ell \neq 3,\ell =1}^5 |K^n_\ell (t,x,y)|
\le C \max \{ t^{-\frac{1}{3}}, t^{-\frac{1}{2}}\}  $ by the above
claims and by Lemma 2.1. Finally we consider $K^n_3.$ Set $K_3=\sum
_n K _3^n $. \proclaim{Lemma  4.6 } There is a fixed $C$ such that
$|   K _3 (t,x,y)| \le C
  \max \{ t^{-\frac{1}{2}}, t^{-\frac{1}{3}} \}  .$
  \endproclaim
{\it Proof.} Set $\chi _{int} (k)=\sum _n \chi _1 (\frac{w - a_{n}^+
}{  |g_n|^{ \frac{1}{4}}}  )  \chi _1 (\frac{ a_{n+1}^--w }{ |g
_{n+1}|^{ \frac{3}{5}}}  )$. Then
$$K _3
(t,x,y)=\int   _{\Bbb R } e^{i(t E  (k) -(x-y)k) }  \chi _{int} (k)
{m ^0_-(x,k)}{ m   ^0_+(y,k)} dk .$$

 We will use the following lemma:

\proclaim{Lemma  4.7}  In the support of $\chi _{int} $ we have
$w\approx k$, $\dot w= 1+O(k^{-2})$,  $\ddot w=  O(k^{-3})$.
Furthermore, we can extend $w$ from the support of $\chi _{int} $ to
the whole of $\Bbb R$ so that the extension (which we denote again
with $w$) satisfies the same relations
\endproclaim
  $w\approx k$ by Lemma 5.2, $\dot w= 1+O(k^{-2})$ by  Lemma 7.1,  $\ddot w=
  O(k^{-3})$ by Lemma 7.4 and the last statement is  Lemma 7.6. Now
  assume Lemma 4.7.
 We have $E=w^2$, $\dot E =2w\dot w$ and $\ddot E=2\dot w^2-2w\ddot
 w$. So by Lemma 4.7 we can extend $E(k)$ from the support of $\chi _{int}$
into a function defined   on all $\Bbb R$       convex with $\dot
E\approx k$, $\ddot E \approx 1$. We express
$$  m _- ^0(x,k)  m  _+ ^0(y,k)  =  \left ( m _- ^0(x,k)  m  _+
^0(y,k)-1  \right )  +1.\tag 1$$

 We have
 Let $k_0$ be the unique
solution of $\dot \Phi (k) =0$. Then set $q^2/2=\Phi (k) -\Phi
(k_0)= \frac{1}{2} \ddot \Phi (\tilde k ) (k-k_0)^2$. Since $1/C \le
\ddot \Phi \le C$ then
  $ \frac{1}{\sqrt{C}}   \le \frac{q}{k-k_0} \le \sqrt{C} $. From
  $q\dot q=\dot \Phi
(k)= \ddot \Phi (k_1) (k-k_0)$ we conclude that for some $C$ we have
$1/C\le \dot q \le C$.  Now we insert (1) in the definition of $K
_3$ obtaining $ K _3 =H_1+H_2  $ with
$$\aligned &
H_1 (t,x,y)=e^{i (tE(k_0)-(x-y) k_0) }
  \int   _{\Bbb R}
e^{itq^2 }  \chi  _{int}(k)\left ( m _- ^0(x,k)  m  _+ ^0(y,k)-1
\right ) \frac{dk}{dq} dq\\& H_2 (t,x,y)= e^{i (tE(k_0)-(x-y) k_0) }
  \int   _{\Bbb R}
e^{itq^2 }
   \chi  _{int}(k)\frac{dk}{dq} dq.
\endaligned
$$
By $\| \widehat{\rho} \| _1\le \tilde C_\varepsilon   \| \rho \|
_{H^{ \frac{1}{2} + \varepsilon}}$,
 in what follows we can use $$ \big | \int   _{\Bbb R} e^{itq^2 }\rho (q)
dq\big | \le C_\varepsilon t ^{-\frac{1}{2}} \|  \rho \| _{H^{
\frac{1}{2} + \varepsilon}}. \tag 4.1$$

\proclaim{Lemma 4.8  }  For a fixed $C$ we have   $\big |  H_2
(t,x,y) \big |  \le  C  t ^{-\frac{1}{2}}$.
\endproclaim
 {  \it Proof . } We write
 $\chi  _{int}(k)= 1-( 1- \chi  _{int}(k)).$ Let $\zeta (k)$ be
 either 1 or $-( 1- \chi  _{int}(k)).$
For $\chi (t)$ a cutoff supported near $t=0$, we
 insert the partition of unity  $\chi (k-k_0)+(1-\chi (k-k_0))$inside
 $ \int    e^{i (tE(k
)-(x-y) k ) }\zeta (k) dk$. Then by $\| \zeta (k)\| _{ \infty  } +
\| \zeta ^\prime (k)\| _{1 }\lesssim 1$ and by Lemmas 3.5 and 4.2,
for a fixed $C$
$$\big | \int e^{i (tE(k )-(x-y) k ) } \zeta (k) \chi (k-k_0)
dk\big | \le C  t ^{-\frac{1}{2}}.$$ Next we want to bound
$$\big | \int e^{i  t\Phi (k) }\zeta (k) (1- \chi (k-k_0))
dk\big |  =\big | \int e^{i \frac{t}{2}q^2 } \zeta (k) (1- \chi
(k-k_0)) \frac{dk}{dq} dq\big | .\tag 2 $$ Let first $\zeta
(k)\equiv 1.$ Then $   \int e^{i \frac{t}{2}q^2 } \zeta (k) (1- \chi
(k-k_0)) \frac{dk}{dq} dq =$ $$=  \int e^{i \frac{t}{2}q^2 } (1-
\chi (k-k_0))   dq  -\int e^{i \frac{t}{2}q^2 }   (1- \chi (k-k_0))
\left (1- \frac{dk}{dq} \right ) dq.$$ By standard arguments the
first term in $O(t^{-\frac{1}{2}})$.  By (4.1), with $\varepsilon
=1/2$, the second term is also $O(t^{-\frac{1}{2}})$ thanks to the
following lemma:

\proclaim{Lemma 4.9 } $\forall \, a>0$ $\exists$ a fixed constant
$C_a$ such that $\| 1- \frac{dk}{dq}\| _{H^1 (\{ |q|\ge a \} }\le
C_a$.
\endproclaim
{\it Proof.} We set $\dot q=\frac{dq}{dk}$. We have $2q\dot q=\dot
E(k)- \dot E (k_0) =2(k-k_0) + O(\langle k \rangle ^{-1}) +
O(\langle k _0\rangle ^{-1}) $ by Lemma 4.1. Hence after integration
$$ \frac{q^2}{(k-k_0)^2}= 1+ \frac{O(\langle k _0\rangle ^{-1})
}{k-k_0 } +\frac{O( \log ( \langle k/k_0 \rangle ))}{ (k-k_0)^2} .
\tag 3$$ We know from $q\approx k-k_0$ that (3) is uniformly bounded
for $|k-k_0|\le 1$. Hence we conclude
$$ 1- \frac{q }{ k-k_0 }=\left (  1+ \frac{q }{ k-k_0 } \right )
^{-1} \left (  1- \frac{q ^2}{ (k-k_0)^2} \right ) $$ is in
$L^2(\Bbb R)$ with norm independent from $k_0$. By $ q\dot q=k-k_0+
O(\langle k \rangle ^{-1})+O(\langle k _0\rangle ^{-1}) $,
$$\dot q =1 +\left (  1- \frac{q }{ k-k_0 } \right ) \dot q + \frac{O(\langle k \rangle ^{-1}) +  O(\langle k
_0\rangle ^{-1})}{k-k_0}.
$$
Since we know $\dot q $ is uniformly bounded, we conclude that for a
fixed constant $C$ we have $\| 1- \dot q \| _{2}\le C$. Next, since
in $\ddot q= \frac{\ddot E -\dot q^{2}}{q}  $ the numerator is
bounded, we see that $\ddot q \in L^2(\{ |q|\ge 1 \} ) .$  By
$\frac{d^2k}{dq^2}= - (\dot q) ^{-3} \ddot q$ we obtain the desired
result.

\bigskip To complete the proof of Lemma 4.8, we have to show that (2) is
$O(t^{-\frac{1}{2}})$ when $ \zeta (k)=-( 1- \chi  _{int}(k)).$ We
have $\| \zeta \| _{1}\lesssim 1$, so for $t\lesssim 1$ we have
$(2)=O(1)$. We suppose now $t\gg 1$. In (2) split $\int _\Bbb R=
\int _{|k|\le \sqrt{t}}+\int _{|k|\ge \sqrt{t}}$. By Lemma 2.1, for
any $N$   $\int _{|k|\ge \sqrt{t}}=O( \sum _{|n|\gtrsim \sqrt{t}}
 |g_n|^{\frac{1}{4}}) \ll t^{-N}.$ We have
$$\int _{|k|\le \sqrt{t}}e^{i  t\Phi (k) }\zeta (k) (1- \chi (k-k_0))
dk=O(    t^{-1}    \sum _{|n|\lesssim \sqrt{t}} (1+1))\approx
t^{-\frac{1}{2}},$$ where the terms 1 are of the form $
 |g_n|^{-\frac{1}{4}}   \int _{a^{+} _{n}} ^ {a^{+}
_{n}+ |g_n|^{ \frac{1}{4}} } dk $ with $ |g_n|^{-\frac{1}{4}}\gtrsim
|\zeta ^\prime | $ in $[a^{+} _{n},  a^{+} _{n}+ |g_n|^{
\frac{1}{4}}]$ and   $ |g_{n+1}|^{-\frac{3}{5}} \int ^{a^{-} _{n+1}}
_ {a^{-} _{n+1}- |g_{n+1}|^{ \frac{3}{5}} } dk $ with $
|g_{n+1}|^{-\frac{3}{5}}\gtrsim |\zeta ^\prime | $ in $[a^{-}
_{n+1}-|g_{n+1}|^{ \frac{3}{5}}, a^{-} _{n+1}]$.

\bigskip

\proclaim{Lemma   4.10 } There is a fixed $C$ such that $|   H_1
(t,x,y)| \le C
  \max \{ t^{-\frac{1}{2}}, t^{-\frac{1}{3}} \}  .$
  \endproclaim
\noindent {\it Proof}. We split $ H_1 (t,x,y)= H _{11} (t,x,y)+H
_{12} (t,x,y)$ with $H _{11} (t,x,y) $ defined inserting   the
additional factor $\widetilde{\chi } _{int}(k)=\sum _n\chi _1 \left
(\frac{k - \pi n }{ {n}^{ \varepsilon -\frac{1}{2}} } \right ) \chi
_1 \left (\frac{(n+1) \pi -k }{{(n+1)}^{ \varepsilon -\frac{1}{2}}}
\right )$ in the definition $H _{1 } (t,x,y) $, and with  $H _{12}
(t,x,y) $ defined inserting $\widetilde{\chi } _{edge}=1-
\widetilde{\chi } _{int}$. The integrals defining $H _{11} (t,x,y) $
are supported in $\pi n +c{n}^{ \varepsilon -\frac{1}{2}}\le  k \le
(n+1) \pi -c(n+1)^{ \varepsilon -\frac{1}{2}}$, that is in the
interior of the bands, while the integrals defining $H _{12} (t,x,y)
$ are supported near the edges.

 \proclaim{Claim} There is a fixed $C$ such that $|H _{11} (t,x,y) |\le C t^{-\frac{1}{2}}$.
 \endproclaim
{\it Proof.} We   have $ {\chi }_{int} \widetilde{\chi
}_{int}=\widetilde{\chi }_{int}$ since $\widetilde{\chi }_{int}$ is
the characteristic function of the  union of $\pi n +c{n}^{
\varepsilon -\frac{1}{2}}\le  k \le (n+1) \pi -c(n+1)^{ \varepsilon
-\frac{1}{2}}$ smoothed and $ {\chi }_{int}$ is the characteristic
function of the  union of  $\pi n + |g_n|^{\frac{1}{4}} \le  k \le
(n+1) \pi -  |g _{n+1}|^{\frac{3}{5}}$ smoothed. We split $H _{11}
(t,x,y) $ in two pieces. For the first piece  we have for $\chi (t)$
a cutoff supported near 0, by Lemma 3.5, by $\| \widetilde{\chi
}_{int} \|_\infty + \| \widetilde{\chi }_{int} ^\prime \chi
(k-k_0)\|_1\lesssim 1$, $ m _- ^0(x,k)  m _+ ^0(y,k)-1= O( k^{-1})$
and $\widetilde{\chi } _{int} (k) \partial _k \left ( m _- ^0(x,k) m
_+ ^0(y,k)\right ) = O( k^{-\frac{1}{2} -\varepsilon}) $,
$$\big |  \int e^{i t\Phi (k) }\widetilde{\chi } _{int} (k)  \chi (k-k_0)
 \left ( m _- ^0(x,k)  m  _+ ^0(y,k)-1  \right )
dk\big | \le C t^{-\frac{1}{2}} .$$ Next we consider

$$  \int e^{i tq^2/2 } \widetilde{\chi } _{int} (k)  (1-\chi (k-k_0)) \left ( m _- ^0(x,k)  m  _+ ^0(y,k)-1  \right )
\frac{dk}{dq}dq.\tag 1$$ By Lemmas 4.5 and 4.9
$$\big \|   \widetilde{\chi } _{int} (k)(1-\chi
(k-k_0))\left ( m _- ^0(x,k)  m  _+ ^0(y,k)-1  \right )
\frac{dk}{dq} \big \| _{H^{1}} \le C    \|  \{   {\langle n\rangle
 ^{ -\frac{1}{2}-\varepsilon}}  \}   \| _{
l^2(\Bbb N )} .
$$
We can apply (4.1) and bound (1) by $C t^{-\frac{1}{2}}$.

\bigskip

\noindent We consider $H _{12} (t,x,y)=\sum _n H _{12} ^{n}(t,x,y)$,
$ H _{12}^{n} (t,x,y)
   =\int _{n\pi}^{(n+1) \pi }
e^{it \Phi ( k)} f(k) dk$

  $$\aligned &\text{with } f(k)=\Psi _n(k)    \left ( m _-
^0(x,k) m _+ ^0(y,k)-1 \right ) \text{ where}\\&  \Psi _n(k)= \chi
_1 (\frac{w - a_{n}^+ }{ |g_n|^{ \frac{1}{4}}} ) \chi _1 (\frac{
a_{n+1}^--w }{ |g _{n+1}|^{ \frac{3}{5}}} )\chi _0 \left (\frac{k -
\pi n }{ {n}^{   \varepsilon -\frac{1}{2}} } \right ) \chi _0 \left
(\frac{(n+1) \pi -k }{{(n+1)}^{ \varepsilon -\frac{1}{2}}} \right
)\endaligned
$$
and  with $  \Phi ( k)= E(k)-E(k_0) - t^{-1} (x-y) (k-k_0)$. Observe
that $ \Psi _n(k)=  \Psi _{n 1}(k) +\Psi _{n 2}(k)$  with $\Psi _{n
1}(k)$ supported in $ |g_n|^{ \frac{1}{4}}\lesssim k-\pi n\lesssim
n^{   \varepsilon -\frac{1}{2}}$ and with $\Psi _{n 2}(k)$ supported
in $ |g _{n+1}|^{ \frac{3}{5 }} \gtrsim  \pi (n+1)-k\gtrsim (n+1)^{
\varepsilon -\frac{1}{2}}.$ Correspondingly write $f=f_1+f_2$ and $H
_{12}^{n}= H _{12}^{n1}+
   H _{12}^{n2} .$
   \proclaim{Lemma 4.11} For a fixed $C$ and  for  $j=1,2$: $|H
   _{12 }^{nj}(t,x,y)|\le C   \langle t  \rangle  ^{-\frac{1}{2}} |\log   t
    |^2 .$
   \endproclaim
{\it Proof.} We  focus on $H _{12}^{n1}$, the proof for  $H
_{12}^{n2}$ being almost the same. We have
$$ H _{12}^{n1} (t,x,y)
   =\int _{n\pi}^{(n+\frac{1}{2}) \pi } F^\prime (k) f_1(k)dk \text{ with } F(k)= \int _{n\pi}^{k}e^{it \Phi ( k^\prime
   )} dk^\prime .$$
   For $H
_{12}^{n2}$ the proof is the same but with  $F(k)= \int _{(n+1)
\pi}^{k}e^{it \Phi ( k^\prime
   )} dk^\prime $.  We get
   $$\aligned & H _{12}^{n1} (t,x,y)
   =-H _{121}^{n1} (t,x,y)-H _{122}^{n1} (t,x,y) \text{ with}\\&
H _{121}^{n1} (t,x,y)=\int _{n\pi}^{(n+\frac{1}{2}) \pi } \left (
\int _{n\pi}^{k}e^{it \Phi ( k^\prime
   )} dk^\prime  \right ) \Psi  _{n1} ^\prime
(k)    \left ( m _- ^0(x,k) m _+ ^0(y,k)-1 \right ) dk\\& H
_{122}^{n1} (t,x,y)=\int _{n\pi}^{(n+\frac{1}{2}) \pi } \left ( \int
_{n\pi}^{k}e^{it \Phi ( k^\prime
   )} dk^\prime  \right ) \Psi  _{n1}   (k)    \partial _k\left
( m _- ^0(x,k) m _+ ^0(y,k)-1 \right ) dk .\endaligned
 $$

\proclaim{Claim} For $ | k_0-n\pi |\ge 2 \pi $ for a fixed $C$  we
have $$|H _{121}^{n1} (t,x,y) |\le C \langle t(\pi n-k_0)  \rangle
^{ -1} \langle  n \rangle ^{ -1}.$$
\endproclaim
{\it Proof.} Indeed we have $|F(k)|\le C \langle t(\pi n-k_0)
\rangle ^{ -1} $, $|m _- ^0(x,k) m _+ ^0(y,k)-1| \le C \langle k
\rangle ^{ -1}$ and $\| \Psi  _{n1} ^\prime (k) \| _{L^1(\pi n , \pi
(n+1))} \le C$. In a similar fashion we obtain \proclaim{Claim} For
$ | k_0-n\pi |\le 2 \pi $ for a fixed $C$  we have $$|H _{121}^{n1}
(t,x,y) |\le C \langle t   \rangle ^{ -\frac{1}{2}} \langle  n
\rangle ^{ -1}.$$
\endproclaim
Next we use that there is a fixed $C$ such that for any $x_0$ and
any $t>0$,
$$\int _{|x-x_0|\ge 1} \frac{dx}{\langle x \rangle \langle t(x-x_0)
\rangle    }  \le  C \min \left \{   t^{-1} , |\log t|\right \} $$
to conclude that for a fixed $C$
$$\sum _n |H _{121}^{n1} (t,x,y) |\le C
\min \left \{   \langle t   \rangle ^{ -\frac{1}{2}}  , |\log
t|\right \} .$$

\bigskip
We now consider $H _{122}^{n1} (t,x,y)$. We start by assuming $ |
k_0-n\pi |\ge 2 \pi $. Then notice that for a fixed $C$
$$\big | \int
_{n\pi}^{k}e^{it \Phi ( k^\prime
   )} dk^\prime \big |  \le \min \{ C \langle t(\pi n-k_0)  \rangle
^{ -1} ,   |k- \pi n|  \} .$$ Next we split $$H _{122}^{n}
(t,x,y)=\int _{n\pi}^{n \pi +\langle t(\pi n-k_0)  \rangle ^{ -1} }
\dots +\int ^{(n+\frac{1}{2})\pi}_{n \pi +\langle t(\pi n-k_0)
\rangle ^{ -1} } \dots
$$
But now
$$ \big |\int _{n\pi}^{n \pi +\langle t(\pi n-k_0)  \rangle ^{ -1} }
\dots  \big | \le C \int _{n\pi}^{n \pi +\langle t(\pi n-k_0)\rangle
^{ -1} }    |k- \pi n| |k- \pi n| ^{-1} \langle  n \rangle ^{ -1}
$$ and
$$ \big |\int ^{(n+\frac{1}{2})\pi}_{n \pi +\langle t(\pi n-k_0)  \rangle ^{ -1}
} \dots  \big | \le   \int ^{(n+\frac{1}{2})\pi}_{n \pi +\langle
t(\pi n-k_0) \rangle ^{ -1} }  C \langle t(\pi n-k_0)\rangle ^{
-1}|k- \pi n| ^{-1} \langle  n \rangle ^{ -1}
$$
and so
$$|H _{122}^{n} (t,x,y)| \le C\langle  n \rangle ^{ -1}\langle t(\pi n-k_0)\rangle ^{ -1} \log \left ( \langle t(\pi n-k_0)\rangle   \right ) .$$
Similarly for $ | k_0-n\pi |\le 2 \pi $ we get
$$|H _{122}^{n} (t,x,y)| \le C\langle t \rangle ^{ -\frac{1}{2}} \log \left ( \langle t \rangle   \right ) .$$
We use that there is a fixed $C$ such that for any $x_0$ and any
$t>0$,
$$\int _{|x-x_0|\ge 1} \frac{\log \left (   \langle t(x-x_0)
\rangle \right ) dx}{\langle x \rangle \langle t(x-x_0) \rangle    }
\le  C \min \left \{  t^{-1} |\log t| , |\log t|^2 \right \}$$ to
conclude that for a fixed $C$
$$\sum _n |H _{122}^{n1} (t,x,y) |\le C \langle t   \rangle ^{
-\frac{1}{2}} |\log t| ^2.$$
\bigskip

\head \S 5 Asymptotic expansion for  $ w-k$
 \endhead
\noindent We consider   $w=u+iv$ and $k=p+iq$. We set $Q_\ell =\frac
1\pi \int _\Bbb R u^\ell q(u) \,du .$ Then $Q_{2\ell +1}=0$. In
particular, see \cite{KK} p. 601, we have $Q_0=\frac 12 \int _0^1
P(t) \, dt $ and $ Q_2=\frac 1 8 \int _0^1 P ^2(t) \, dt.$

For $u\in   \sigma  $ we have $q(u)=0.$ For $u\in g_n$ we have
formula (4.12) \cite{K1}:
$$q(u)=\sqrt{ (u-a_n^-) (  a_n^+-u  )}  \left ( 1+\frac 1\pi  \sum _{m\neq n}\int
_{g_m} \frac{q(t) \, dt }{|t-u| \sqrt{ (t-a_n^-) (  t-a_n^+    )}}
\right ) .\tag 5.1$$ By Lemma 2.1 there is a fixed constant $C>0$
such that for $|\sigma _m|$ the length of $ \sigma _m $ we have $
|\sigma _m| \ge C $ for all $ m.$ For $u\in g_n$ and $ C_0=
 1+    \frac{Q_0}{\min \{ |\sigma _m| \colon m\in \Bbb Z \} }$ by
 \cite{K1} p.16
$$\aligned &     \sqrt{ (u-a_n^-) (  a_n^+-u  )} \le
  q(u)\le  C_0\sqrt{ (u-a_n^-) (  a_n^+-u  )}.
 \endaligned \tag 5.2$$

\proclaim{Lemma 5.1} 1) $\forall$  $N$ there is a $C_N>0$ such that
$|q(u)| \le C_N \langle u \rangle ^{-N}$ $\forall$  $u\in \Bbb R$.

\noindent 2) The distributional derivative $q^\prime (u)$ satisfies
$\langle u \rangle ^Nq^\prime \in L^{ r}(\Bbb R)$ for any $1\le r<2
$ and any $N$.
\endproclaim
{\it{Proof.}} We have $ 0\le q(u)\le C_0|g_n| \le C_N \langle u
\rangle ^{-N}$ for $u\in g_n$  by (5.2) and Lemma   2.1.    Turning
to the second claim, by (5.1) the pointwise derivative $q^\prime
(u)$ is well defined except at the points $a_n^\pm $ for $n\in \Bbb
Z$. Obviously $q^\prime (u)=0$ for any $a_n^+<u<a_{ n+1} ^{-}$. For
$a_n^-<u<a_{ n } ^{+}$ we differentiate   (5.1) and using the fact
that inside the integral we have $|t-u| \ge \inf _n|\sigma _n| >0$,
we conclude there is a fixed $C$ such that
$$|q^\prime (u)|< C\left ( \frac{a_{ n } ^{+}-u}{u-a_{ n }
^{-}}\right ) ^{\frac{1 }{2}} + C \left ( \frac {u-a_{ n } ^{-}}
{a_{ n } ^{+}-u}\right ) ^{\frac{1}{2}} .
$$
From this we conclude that the pointwise $q^\prime (u)$ coincides
with the   distributional derivative and that $\int _{ g_n}\langle u
\rangle ^N|q^\prime (u)|^rdu\le \frac{C\langle n \rangle ^N}{2-r}
|g_n|$ for some fixed $C$. By Lemma 2.1 we conclude $\| \langle u
\rangle ^Nq^\prime \| _r \le \frac{C_N}{2-r}$ for some $C_N$.

\bigskip

From the second claim   in Lemma 5.1  we obtain:

 \proclaim{Lemma  5.2} For any integer $N$
there is a constant $C_N>0$ such that, for any $w=u+iv$ with   $
v\ge 0$ and $|u|>1$, we have
$$w-k(w) = \sum _{\ell =0} ^{N}  \frac{ Q_\ell  }{w^{\ell +1}}+
R_N(w)\, , \quad \big | R_N(w) \big | \le   \frac{C _N}{\langle w
\rangle ^{N+1} }.$$
\endproclaim

{\it{Proof}}.  We have   for $v>0$ by (4.1) \cite{K1} $k(w)-w= \frac
1\pi \int _{\cup g_n} \frac{ q(t)   }{t-w }       dt$ and so
$$ w-k(w)=  \frac 1\pi \sum _{\ell =0} ^N \frac{1}{
w ^{\ell+1}}
        \int _{\cup g_n} t^\ell  q(t)
 dt + \frac 1{  w^{N+1}}  \left (  Q_v - i   P_v \right ) \ast (t^{N+1 }   q(t) ) (u)
   $$
with in the last term  the convolution of $t^{N+1 }   q(t)$ with the
Poisson kernels
$$P_v(x)= \frac{1}{\pi}\frac{v}{x^2+v^2} \, , \quad  Q_v(x)=\frac{1}{\pi} \frac{x}{x^2+v^2}.$$
Since   for $1<r<2$ there is a $C_r$ such that for any $v> 0$, see
p. 121 \cite{Ste},
$$\| \left (  Q_v - i   P_v \right ) \ast (t^{N+1 }   q(t) ) \| _{W^{1,
r}} \le C_r\| t^{N+1 }   q(t)  \| _{W^{1, r}}\, ,
$$
by Lemma 5.1 and by the  Sobolev embedding theorem there is  a fixed
constant $C_N>0$ such that  $|\left (  Q_v - i   P_v \right ) \ast
(t^{N+1 } q(t) ) (u)|\le C_N. $

\bigskip
As an immediate corollary of Lemma 5.2 we obtain:

\proclaim{Lemma 5.3} There are two constants $C_1$ and $C_2$ such
that, for $w=u+iv$ with $1\ge v\ge 0$ and $|u|>C_1$, we have
$$\big | u-p-   \frac{ Q_0  }{u}+ \frac{ Q_0  v^2}{u^3}-\frac{ Q_2  }{u^3}\big |  \le   \frac{C _2}{|u| ^{4} }.$$
\endproclaim
\noindent By $w-k= Q_0 \frac{  \overline{w}}{|w|^2} + Q_2\frac{
\overline{w}^3}{|w|^6}+ O(w^{-4})$, we get $u-p=Q_0 \frac{
u}{u^2+v^2} + Q_2\frac{ u^3}{(u^2+v^2)^3}+ O(u^{-4})$. Then by
Taylor series $\frac{ 1}{u^2+v^2} = u^{-2}- v^2 u^{-4}+\dots ,$ we
get the desired result.

\bigskip \head \S 6 Relation between $v$ and $q(u+iv)$ \endhead

In the proof of Lemma 4.4  we will need to use the  relative size of
the   coordinates of $w=u+iv$ and $k(w)=p(w)+iq(w).$ Lemma 5.3 gives
some information on $u-p(u+iv)$. We now consider the relation
between $v$ and $q(u+iv)$. Recall that
$$q(u+iv)=v+\frac{v}{\pi}\int _\Bbb R\frac{q(t)\, dt}{(t-u)^2+v^2}.
\tag 6.1
$$
 The unnatural restriction on $u$ in the following lemma is not
sharp and is only justified by estimates needed later.

 \proclaim{Lemma 6.1} Let $a^{+}_{n}+  |g_n|^{\frac{1}{4}}\le u
\le a^{-}_{n+1}-  |g _{n+1}|^{\frac{3}{5}}$ for   $n \gg 1$ large,
and let $1\ge v \ge 0$ Then there is a constant $C $   such that
$$|q(u+iv)-v-\frac{v}{u^2}Q_0|\le C v u^{-4}.$$

\endproclaim

{\it Proof.} By (6.1)  write
$$q=v+\frac{v }{\pi u^2} \int _\Bbb R q(t) \left ( 1+\frac{2t}{u}-
\frac{t^2}{u^2}-\frac{v^2}{u^2} \right ) dt +\frac{v }{\pi  } \int
_\Bbb R \frac{q(t) \left (  \frac{2t}{u}-
\frac{t^2}{u^2}-\frac{v^2}{u^2} \right )^2}{(t-u)^2+v^2} dt
$$
and use the formulas for $Q_\ell $ above (5.1) to express the second
term on the rhs   as $\frac{v }{  u^2} (
(1-\frac{v^2}{u^2})Q_0-\frac{Q_2 }{u^2}).$ In  the second integral
we expand the square on the numerator, treating the resulting terms
separately. For example we write
$$\aligned   \frac{4v }{\pi  u^2} \int
_\Bbb R \frac{t^2q(t)  }{(t-u)^2+v^2} dt& = \frac{4v }{\pi  u^4}
\int _\Bbb R  t^2q(t)  \left ( 1+\frac{2t}{u}-
\frac{t^2}{u^2}-\frac{v^2}{u^2} \right )   dt+\\& +
 \frac{4v }{\pi  u^4} \int
_\Bbb R \frac{t^2q(t)  \left (   2t  - \frac{t^2}{u }-\frac{v^2}{u }
\right )^2 }{(t-u)^2+v^2} dt.\endaligned $$ The first term on the
rhs is $O(vu^{-4})$.  To show that the second term is $O(vu^{-4})$,
we need to show bounds of the form $ \int _\Bbb R \frac{|t|^N q(t)
}{(t-u)^2+v^2} dt \le C_N.$ Say that $u\le
\frac{a^{+}_{n}+a^{-}_{n+1}}{2}$, with the case $u\ge
\frac{a^{+}_{n}+a^{-}_{n+1}}{2}$ treated similarly. Then
$$ \int _\Bbb R \frac{|t|^N q(t)
}{(t-u)^2+v^2} dt= \left ( \sum _{\ell \neq n }  \int _{g_\ell
}\frac{|t|^N q(t) }{(t-u)^2+v^2} dt \right ) + \int _{a^{-}_{n} }
^{a^{+}_{n} }\frac{|t|^N q(t) }{(t-u)^2+v^2} dt  $$ where the first
term in the rhs is bounded by a $C_N\langle u \rangle ^{-2}$ thanks
to $|t-u|\approx |\ell -n |$, Lemma 5.1 and $ \langle u \rangle
 ^{-2}\ast \langle u \rangle  ^{-N}\lesssim \langle u \rangle  ^{-2}.$ Next we
write
$$\int _{a^{-}_{n} }
^{a^{+}_{n} }\frac{|t|^N q(t) }{(t-u)^2+v^2} dt  \le \frac{C|n|^N
|g_n|^2}{(a^{+}_{n}-u)^2}\le C_N
$$
by    $q(t)\le (1+C n^{-2}) |g_n|$, see \cite{KK}, by $a^{+}_{n}
-a^{-}_{n}=|g_n|$, by our restriction on $u$ and by Lemma 2.1.

\bigskip
We have the following corollary:
 \proclaim{Lemma 6.2} In   $\{ 1\ge  v \ge 0\} $ for any
preassigned $C_1>0$ in $a^{+}_{n}+ |g_n|^{\frac{1}{4}}\le u \le
a^{-}_{n+1}- |g _{n+1}|^{\frac{3}{5}}$ there is a    $C$ such that

   $$ \big | \Im
(E-k^2)  \big | \le C \frac{  v}{u^3}    .$$

\endproclaim
\noindent {\it Proof.} Write $\Im (E-k^2) =(v-q)(u+p)+(u-p)(v+q) $
with $0\le v \le q$. By Lemma 5.3 we have
$u-p=\frac{Q_0}{u}-\frac{Q_0v^{2}}{u^3}+\frac{Q_2}{u^3} + O(u^{-4} )
$ and $u+p=2u- \frac{Q_0}{u}+\frac{Q_0v^{2}}{u^3} -\frac{Q_2}{u^3} +
O(u^{-4} ). $ By Lemma 6.1 $v-q=-vu^{-2}Q_0+O( vu^{-4})$ and
$q+v=2v+vu^{-2}Q_0+O( vu^{-4})$. The $O(vu^{-1 })$ term in $\Im
(E-k^2)$ cancels and   the following one is $O(vu^{-3 })$.

\bigskip

\noindent Now we extend Lemma  6.2 without the restriction on $u$.
From $q(t)\ge 0$ and $q(t)\not  \equiv 0$ we get $ q(u+iv)> v$ by
(6.1). The following is an elementary consequence of Lemmas 2.1 and
5.1:

\proclaim{Lemma 6.3} For any $u\in \Bbb R$ there is at most one
$g_n$ such that $\text{dist}(u,g_n)\ll 1$. For such an $n$ we have
$$ 0\le \sum _{\ell \neq n}\frac{v}{\pi}\int  _{g_\ell }\frac{q(t)\,
dt}{(t-u)^2+v^2} <C \frac{v}{\langle u \rangle ^2}.$$ If such an $n$
does not exist the above formula holds summing over all $\ell \in
\Bbb Z$.
\endproclaim

\noindent Suppose now that $u$ is close to the gap $g_n=(a^{-} _{n},
a^{+} _{n})$  and   set $$I_n(u,v )=\frac{v}{\pi}\int _{a^{-}
_{n}}^{a^{+} _{n}} \frac{q(t)\, dt}{(t-u)^2+v^2}.  $$

\proclaim{Lemma 6.4} There is a fixed $C$ independent from $n$ such
that:

\roster

\item Let $u\in \sigma _n.$ Set $m=n,n+1$, $a_m$ either $a_m=a_n^+$
for $m=n$ or $a_m=a_{n+1}^-$ for $m=n+1$. Suppose $w$ close to
$a_m$. Then
 if
$|u-a  _{m}|>\frac{|g_m|}{4}$ or if $v\ge |g_m|$    we have
$$\frac 1C \frac{|g_m|^2 v}{(u-a  _{m})^2+v^2}< I_n(u,v )< C \frac{|g_m|^2 v}{(u-a  _{m})^2+v^2} ; $$

\medskip

\item  if $u\in g_n=[a^{-} _{n}, a^{+} _{n}]$  for $v\ge \frac{|g_n|}{2}$ we have

$$ \frac 1C \frac{|g_n|^2 }{v}<I_n(u,v )<C \frac{|g_n|^2 }{ v } ;$$

\medskip

\item for $u\in g_n$ and $ v<2\min
  \{   |u-  a^{-} _{n}|,    |u-  a^{+} _{n}| \}
 $ we have $$\frac 1C |g_n|^{ \frac{1}{2}} \sqrt{\min
  \{   |u-  a^{-} _{n}|,    |u-  a^{+} _{n}| \}}< I_n(u,v )<C
|g_n|^{ \frac{1}{2}} \sqrt{\min   \{   |u-  a^{-} _{n}|,    |u-
a^{+} _{n}| \}}.$$

\endroster

\endproclaim
 {\it Proof}. {\item{(1)}} For definiteness let $u\le
 (a_n^++a_{n+1}^-)/2$.
Set $\tilde u = u- a^{+} _{n} $.   We have
$$I_n(u,v)= \frac{v}{\pi}\left [ \int  _{0}^{ \frac{|g
_{n}|}{2}} +  \int  _ { \frac{|g _{n}|}{2}}  ^{|g _{n}|} \right ]
\frac{ \sqrt{t} \, \sqrt{  |g _{n}|-t } }{(t+\tilde u)^2+v^2} dt.
$$
We have
$$ v \int  _ { \frac{|g _{n}|}{2}}  ^{|g _{n}|}
\frac{ \sqrt{t} \, \sqrt{  |g _{n}|-t } }{(t+\tilde u)^2+v^2}
dt\approx \frac{|g _{n}|^2 v}{(|g _{n}|+\tilde u)^2+v^2}
$$
Set  $\tilde I_n(u,v)=v \int  ^ { \frac{|g _{n}|}{2}}  _0 \frac{
\sqrt{t} \, \sqrt{  |g _{n}|-t } }{(t+\tilde u)^2+v^2} dt$. For
$\tilde u \gtrsim |g _{n}|$ or for $ v\gtrsim |g _{n}|$ then
$$\tilde I_n(u,v)=\frac{v}{ \tilde u ^2+v^2}\int  ^ { \frac{|g _{n}|}{2}}  _0\sqrt{t} \, \sqrt{  |g _{n}|-t
}dt\approx \frac{|g _{n}|^2 v}{ \tilde u ^2+v^2} .
$$

 {\item {(2)}} For definiteness here and below $u\le \frac{a^{+} _{n} +a^{-}
_{n}}{2}$. We  set $\tilde u = u- a^{-} _{n} $.   We have by (5.2)
$$I_n(u,v)\approx \frac{v}{\pi}\int  _{a^{-} _{n}}^{a^{+}
_{n}} \frac{ \sqrt{a^{+}_n-t} \, \sqrt{t- a^{-}_n } }{(t-u)^2+v^2}
dt\approx v \sqrt{|g _{n}| }\int  _{0}^{\frac{|g _{n}|}{2}} \frac{
   \sqrt{t } dt}{(t-\tilde u)^2+v^2}
$$
Then for $v\gtrsim |g_n|$  we get $I_n(u,v)\approx \frac{|g_n|^2v}{
\tilde u ^2+v^2}$ and hence claim 1 in Lemma 6.4.

{\item{(3)}} We write for $\tilde u = u- a^{-} _{n} $

$$\aligned &I_n(u,v)\approx v \sqrt{|g _{n}| } \left [ \int  _{0}^{\frac{\tilde
u}{2}   } + \int  _ {\frac{\tilde u}{2}   }^{2\tilde u}  +\int
_{2\tilde u }^{\frac{|g _{n}|}{2}} \right ] \frac{ \sqrt{t }
dt}{(t-\tilde u)^2+v^2} \approx \frac{v \sqrt{|g _{n}|} \tilde u^{
\frac{3}{2}}} { \tilde u ^2+v^2} +\\& +  \sqrt{|g _{n}| \tilde u}
\arctan (\frac{\tilde u}{2v})  +   \sqrt{|g _{n}|v}  \int
_{\frac{2\tilde  u }{v}}^{\frac{|g _{n}|}{2v}}  \frac{ \sqrt{t }
dt}{t^2+1}  .
\endaligned
$$
For $2\tilde  u> v$ we get $I_n(u,v)\approx \sqrt{|g _{n}|\tilde u}$
and hence claim 3 in Lemma 6.4.

\bigskip \head \S 7 Estimates on the band function \endhead

We will need to bound $\frac{d E}{dk }$,  $\frac{d^2E}{dk^2}$ and
$\frac{d^3E}{dk^3}$.

\proclaim{Lemma  7.1} There are  constants $C_1>C_2>0$ such $\forall
m$ and    $\forall u\in \sigma _m=[a^+_m,a^-_{m+1}]$ and if $v=0$,
we have for $A(u)= \frac{ |g_m|^2   }{(
 u-a_m^+    )^{\frac{1}{2}}  ( u-a_m^++|g_m|)^{\frac{3}{2}}  }  +  \frac{ |g_{m+1}|^2
}{(
  a_{m+1}^- -u   )^{\frac{1}{2}}  (  a_{m+1}^--u+|g _{m+1}|)^{\frac{3}{2}}  }$
  $$ \aligned &  1+  C_2\left ( A(u) + \frac{1}{\langle u \rangle ^2} \right ) \ge    p^\prime (u) \ge  1+C_1 A(u)    .\endaligned \tag 1$$
  Correspondingly for $p\in [m\pi , (m+1) \pi ]$ we have
  $$   \align &\frac 1{1+C_2 \left ( A(u) + \frac{1}{\langle u \rangle ^2 } \right )}\le       \frac{du}{dp}\le       \frac 1{1+C_1 A(u)}
   \tag 2\\& \frac
    {2\big |u \big | }{1+C_2 \left ( A(u) + \frac{1}{\langle u \rangle ^2} \right )} \le \big | \frac{dE}{dp}\big | \le
       \frac{2\big |u \big | }{1+C_1 A(u)}  .\tag 3\endalign$$
        If $a^{+}_{n}+  |n| ^{3} |g_n| \le u
\le a^{-}_{n+1}-  |{n}+1|^{3} |g _{n+1}|$  with $n \gg 1$ large and
$v=0$, then there is a fixed $C$ such that
$$|\dot w-1|+|\dot E (k)-2k | \le C \langle k \rangle ^{-2}. \tag 4$$
There is $n_0$ such that if for $n\ge n_0$, $ |a_{n}^{+}-u| +|v|\ge
|k| |g_n| $ with  $a^{+}_{n} \le u \le
\frac{a^{+}_{n}+a^{-}_{n+1}}{2} $ and $1\ge  v\ge 0 $ then there is
a $C$ such that for the corresponding $k=p+iq$ we have $$ \big |
\frac{\dot E}{2k}-1 \big | \le C |k|^{-1}
 .\tag 5$$
         (5) holds also in $ |a_{n+1}^{-}-u| +|v|\ge |k|
|g_{n+1}| $ and in particular     for $a^{+}_{n}+
|g_n|^{\frac{1}{4}}\le u \le a^{-}_{n+1} -|g_{n+1}|^{\frac{3}{5}} $
and $1\ge  v\ge 0 $.

\endproclaim
Estimate  (5) is used in the proof of Lemma 4.4.

 {\it Proof of Lemma 7.1. } By \cite{K1}, $ p^\prime (u)=1+ \frac
1\pi \sum _n \int _{ g_n}    \frac{ q(t)   }{(t-u)^2}       dt\ge
1+I(u) $ with
$$I(u)=  \frac 12 \sum _n \frac{|g_n|^2 }{\sqrt{
|(u-a_n^-) (  a_n^+-u  )|    }            \left ( \sqrt{ | u-a_n^-
|} +   \sqrt{ |     a_n^+-u  | } \right ) ^2} .$$   For $u\in \sigma
_m$, ignoring all the terms
 in the sum defining $I(u)$
 except for
$n=m,m+1$, we get the lower bound for $p^\prime (u)$ in (1) by $
\sum _{n = m, m+1 }\cdots \approx A(u)$. Turning to the  upper
bound, by (5.2) we have $p^\prime (u)\le 1 + C_0I(u)$. We split now
$2I(u)=$
$$ \left (  \sum _{n \neq m, m+1 } +   \sum _{n = m, m+1 }        \right ) \frac{|g_n|^2 }{\sqrt{
|(u-a_n^-) (  a_n^+-u  )|    }            \left ( \sqrt{ | u-a_n^-
|} +   \sqrt{ |     a_n^+-u  | } \right )^2} .$$ Since there is a
fixed $c>0$ such that for  $n \neq m, m+1$ and for $u\in \sigma _m$
we have $| a_n^\pm -u  |\ge c \langle n-m \rangle $ then
$$\sum _{n \neq m, m+1 }\cdots \le C \sum _{n   }  \frac{|g_n|^2 }{\langle
n-m \rangle ^2}   \le C_2 \frac{1}{\langle u \rangle ^2}  .$$
 This gives the
upper bound for $p^\prime (u)$ in (1). (2) is obtained taking the
inverses in (1) and (3) follows from $\frac{d E}{dk }= 2w\frac{d
w}{dk }$.

To prove (4) we claim \proclaim{Claim} For   $a^{+}_{n}+C_1 |n| ^{3}
|g_n| \le u \le a^{-}_{n+1}-C_1 |{n}+1|^{3} |g _{n+1}|$ for any
fixed $C_1$, $n \gg 1$ large and $v=0$,  we have $ p^\prime (u)=1+
\frac{Q_0}{u^2} +O(u^{-3})$. \endproclaim \noindent Assume for a
moment the Claim. Then $\dot w =1- \frac{Q_0}{u^2}+O(u^{-3})$ and by
Lemma 5.3
$$\dot E= 2u \dot w= 2 u-2\frac{Q_0}{u}+O(u^{-2})=   2p +2\frac{Q_0}{u}- 2\frac{Q_0}{u} +O(u^{-2}).$$

\noindent{\it Proof of the Claim.} Is suggested by formal
differentiation of  $p(u)=u -\frac{Q_0}{u}-\frac{Q_2}{u^3}+\dots $,
but for a proof
 we return to $ p^\prime (u)=1+ \frac 1\pi \sum _n \int _{ g_n}
\frac{ q(t)   }{(t-u)^2}       dt$.  By the argument   in the proof
of Lemma 6.1, simply setting $v=0$ in the appropriate integral,
$$p^\prime (u)=1+ \frac{Q_0}{u^2}- \frac{Q_2}{u^4}+\frac{1}{\pi u^2}\int
_\Bbb R \frac{q(t) \left (   2t - \frac{t^2}{u}  \right )^2}{(t-u)^2
} dt .
$$ Then
$$ \int
_\Bbb R \frac{q(t) \left (   2t - \frac{t^2}{u}  \right
)^2}{(t-u)^2}= \left [ \left ( \sum _{\ell \neq n }  \int _{g_\ell }
\right ) + \int _{a^{-}_{n} } ^{a^{+}_{n} }  \right ] \frac{q(t)
\left ( 2t - \frac{t^2}{u} \right )^2}{(t-u)^2} dt  $$ where the
first term in the rhs is bounded by a $C \langle u \rangle ^{-2}$
thanks to $|t-u|\approx |\ell -n |$ and where
$$\int _{a^{-}_{n} }
^{a^{+}_{n} }\frac{q(t) \left ( 2t - \frac{t^2}{u} \right
)^2}{(t-u)^2} dt  \lesssim \int _{a^{-}_{n} } ^{a^{+}_{n}
}\frac{q(t) t  ^2}{(t-u)^2} dt \le \frac{C |n|^2|g_n|
^2}{(a^{+}_{n}-u)^2} \lesssim   n^{-4}
$$
   by $t\approx u\approx n$ i.e. Lemma 2.1,   $q(t)\le (1+C n^{-2}) |g_n|$, see \cite{KK},  $a^{+}_{n}
-a^{-}_{n}=|g_n|$, see definition of $|g_n|,$ and by our restriction
on $u$, that is $u-a^{+}_{n}\gtrsim n^3|g_n| $.

To prove (5) we write $$ k^\prime (w)=1+ \frac 1\pi \sum  _{\ell
\neq n} \int _{ g_\ell } \frac{ q(t)   }{(t-w)^2}       dt+   \frac
1\pi \int _{a_n^-}^ {a_n^+} \frac{ q(t)   }{(t-w)^2}       dt . $$
The second term in the rhs is $O(n ^{-2})$. For
$|a_{n}^{+}-u|+|a_{n}^{-}-u| +|v|\ge |k| |g_n|, $  so in particular
for  $a^{+}_{n}+c {n} ^{\frac{1}{3}} |g_n|^{\frac{2}{3}}\le u \le
\frac{a^{+}_{n}+a^{-}_{n+1}}{2} $ and $1\ge  v\ge 0 $,   the third
term has absolute value less than
$$ \frac 1\pi \int
_{a_n^-}^ {a_n^+} \frac{ q(t)   }{(t-u)^2+v^2}       dt\le C \frac{
|g_n|^2} {|k|^2     |g_n| ^{2} } = O(k^{-2}).$$ So $\dot w =
1+O(k^{-2})$  and $\dot E=2w\, \dot w=2w+O(k^{-1})= 2k+O(k^{-1})$ by
Lemma 5.2.

\bigskip

We will need the following formulas, see (6.1)
\cite{K1}:\proclaim{Lemma 7.2} For any $u\not \in [a,b]$ we have

$$\align &
\int _{a}  ^b  \frac{ \sqrt{(t-a)(b-t)}   }{(t-u)^3} dt=\frac{\pi
}{8}\text{sign} (u-a) \frac{(b-a)^2}{|u-a|^{\frac 32}|u-b|^{\frac
32}} \tag 1 \\& \int _{a}  ^b  \frac{ \sqrt{(t-a)(b-t)}   }{(t-u)^4}
dt=\frac{ \pi }{16}  \frac{(b-a)^2}{|u-a|^{\frac 32}|u-b|^{\frac
32}}     \left (    \frac{1}{|u-a| }   +   \frac{1}{|u-b| }  \right
)        .  \tag 2
\endalign  $$

\endproclaim
The proof is as Lemma 6.1 \cite{K1}. Let us suppose $u>b$. Then
setting $\ell =b-a$,  $\ell s= t-a$, $\ell h=a-u$, $q = \sqrt{1+
\frac{1}{h }}$ and introducing a new variable $x$ defined by $x\, s
= \sqrt{s(1-s)}$ and so $s=1/(1+x^2),$ $ds=-2xs^2dx$, we obtain
$$\frac{1}{\ell }\int _0^1 \frac{\sqrt{s(1-s)}}{(s+h)^3}ds=\frac{1}{\ell h^3}\int _0^1
\frac{\sqrt{s(1-s)}}{s^3(\frac 1s+\frac 1h)^3}ds =\frac{1}{\ell
h^3}\int _\Bbb R \frac{x^2}{(x^2+q^2)^3}dx$$ which by the Residue
Theorem is equal to $
 \frac{ \pi i }{\ell h^3}  \left [     \frac{z^2}{(z+iq)^3} \right
] ^{\prime\prime} _{z=iq } = \frac{ \pi  }{8\ell h^3 q^3} . $
Proceeding similarly we get (1) also for $u<a$.  (2) follows by
differentiation.

\bigskip

\proclaim{Lemma  7.3} There are positive constants $C_0$, $C_1$,
$C_2$, $C_3$,  $\alpha $ and $m_0\ge 0$ such that    for any $ m\ge
m_0$ and for any $u\in ]a^+_{m }   , a^-_{m+1} [$ we have

$$\align &   |p^{\prime \prime }(u)|\le C_1   \langle u \rangle  ^{-3}\quad \forall u\in  [a^+_m+\alpha ,a^-_{m+1} -\alpha
] \tag 1 \\& p^{\prime \prime }(u)\le   -  C_1\langle u \rangle
^{-3}- \frac{ 1}{4}\frac{(a_m^+-a_m^-)^2}{|u-a_m^-|^{\frac
32}|u-a_m^+|^{\frac 32}} \quad \forall u \le  a^+_m+\alpha  \tag 2
\\&
p^{\prime \prime }(u)\ge   -  C_2\langle u \rangle ^{-3}-
C_0\frac{(a_m^+-a_m^-)^2}{|u-a_m^-|^{\frac 32}|u-a_m^+|^{\frac 32}}
\quad \forall u \le  a^+_m+\alpha  \tag 3
\\&p^{\prime \prime }(u) \ge - \frac{C_2}{\langle u
\rangle ^{ 3}}+\frac {1}4
\frac{(a_{m+1}^+-a_{m+1}^-)^2}{|u-a_{m+1}^-|^{\frac
32}|u-a_{m+1}^+|^{\frac 32}}  \quad \forall u\ge a^-_{m+1} -\alpha
\tag 4  \\& p^{\prime \prime }(u) \le - \frac{C_3}{\langle u \rangle
^{ 3}}+\frac {C_0}4
\frac{(a_{m+1}^+-a_{m+1}^-)^2}{|u-a_{m+1}^-|^{\frac
32}|u-a_{m+1}^+|^{\frac 32}}  \quad \forall u\ge a^-_{m+1} -\alpha
.\tag 5
\endalign$$
Since $p(u)$ is odd, for $m\le - m_0$ there is an analogous
statement.
\endproclaim

We start with
$$p^{\prime \prime }(u)=  \frac 2\pi \sum _n \int _{ g_n}    \frac{ q(t)   }{(t-u)^3}       dt   .$$
We are assuming $u\in \sigma _m=[a_m^+,a_{m+1}^-] $. The terms with
$n\le m$ (resp. $n> m$) are negative (resp. positive). We have

$$\aligned &p^{\prime \prime }(u) \le    \frac 2\pi \sum _{n\le m}
\int _{ g_n} \frac{ \sqrt{  (t-a_n^-) (  a_n^+-t  ) }   }{(t-u)^3}
dt +\frac {2C_0}\pi \sum _{n> m}\int _{ g_n}    \frac{ \sqrt{
(t-a_n^-) ( a_n^+-t  ) }   }{(t-u)^3}       dt \\& = -\frac 14 \sum
_{n\le m}  \frac{(a_n^+-a_n^-)^2}{|u-a_n^-|^{\frac
32}|u-a_n^+|^{\frac 32}} +\frac {C_0}4 \sum _{n> m}
\frac{(a_n^+-a_n^-)^2}{|u-a_n^-|^{\frac 32}|u-a_n^+|^{\frac 32}}.
\endaligned \tag 6
$$

Similarly we have
$$\aligned &p^{\prime \prime }(u) \ge   -\frac {C_0}4
\sum _{n\le m}  \frac{(a_n^+-a_n^-)^2}{|u-a_n^-|^{\frac
32}|u-a_n^+|^{\frac 32}} +\frac {1}4 \sum _{n> m}
\frac{(a_n^+-a_n^-)^2}{|u-a_n^-|^{\frac 32}|u-a_n^+|^{\frac 32}}.
\endaligned \tag 7
$$
We are considering  $m\gg 1$. Observe that by Lemma 2.1 for $u\in
\sigma _m$ there
  are fixed constants   such that for arbitrary $N$

$$\aligned & \sum _{n\ge  m+1}
\frac{(a_n^+-a_n^-)^2}{|u-a_n^-|^{\frac 32}|u-a_n^+|^{\frac
32}}\approx  \sum _{n> m+1}  \frac{|g_n|^2}{\langle n-m \rangle
^{3}}\le \frac{C_N }{\langle  u \rangle ^{N}} \\& \sum _{n\le m-1}
\frac{(a_n^+-a_n^-)^2}{|u-a_n^-|^{\frac 32}|u-a_n^+|^{\frac 32}}
\approx    \sum _{n\le m-1} \frac{|g_n|^2}{\langle n-m \rangle
^{3}}\approx   \langle u \rangle ^{-3}  . \endaligned \tag 8$$
 Hence for $u\in  [a_m^++ \alpha ,a_{m+1}^--\alpha ] $
with $\alpha
>0$,
 (6)-(8) and $|g _{n}|\lesssim \langle n \rangle ^{-N}$
imply (1). Assume now   $u\in (a^+_m , a^+_m+\alpha ]$. From
(6)-(8), $|g _{m+1}|\lesssim \langle m+1 \rangle ^{-N}$ and the
signs, we get (2) and (3)
$$ \aligned &p^{\prime \prime }(u) \approx  -\langle u \rangle ^{-3} -  \frac{(a_m^+-a_m^-)^2}{|u-a_m^-|^{\frac
32}|u-a_m^+|^{\frac 32}}   .
\endaligned
$$

Now we consider $u$ close to $a^{-}_{n+1}$. We now prove (4). From
(7) we get

$$ \aligned &p^{\prime \prime }(u) \ge -\frac {C_0}4 \sum
_{n\le m}  \frac{(a_n^+-a_n^-)^2}{|u-a_n^-|^{\frac
32}|u-a_n^+|^{\frac 32}} +\frac {1}4 \sum _{n> m+1}
\frac{(a_n^+-a_n^-)^2}{|u-a_n^-|^{\frac 32}|u-a_n^+|^{\frac 32}}
+\\&    +\frac {1}4
\frac{(a_{m+1}^+-a_{m+1}^-)^2}{|u-a_{m+1}^-|^{\frac
32}|u-a_{m+1}^+|^{\frac 32}} .
\endaligned
$$
We   absorb the first two terms in the right hand side inside the
term $\frac{-C_2}{\langle u \rangle ^3}$ of (4) and we get (4). The
proof of (5) proceeds similarly starting from (6).

\bigskip \noindent In the following two lemmas the symbols $\approx $,
$\lesssim $ and $\ll $ involve fixed constants. We remark that
$E(k)$ is   even in $k$ so for this reason  we will assume now only
$k\gg 1$.

\proclaim{Lemma 7.4} There are   fixed $C_1,C_2,c>0$, with
$C_1>C_2$,  and $m_0>0$ such that for any $m\ge m_0>0$ and any $u\in
\sigma _m$ we have $ \big |\frac{d^2E}{du^2}\big | \approx $
$$\aligned   \frac{m}{|g_{m  }|} \quad &\text{for}\quad u-a_m^+
 \le  c|g_{m  }|; \\  1+       \frac{m|g_m|^2}{ |u-a_m^+|^3}        \quad &\text{for}\quad
a_m^++c  |g_{m  }|\le u \le a_{m+1 }^--C
_1|m+1|^{\frac{1}{3}}|g_{m+1 }|^{\frac{2}{3}};
\\    \frac{-m|g_{m+1 }|^2}{|u-a_{m+1 }^-|^{\frac
32}|u-a_{m+1 }^+|^{\frac 32}}  \quad &\text{for}\quad  c|g_{m+1
}|\le a_{m+1 }^--u \le C_2  |m+1|^{\frac{1}{3}}|g_{m+1
}|^{\frac{2}{3}};
\\
 \frac{-m}{|g_{m+1 }|}
\quad &\text{for}\quad a_{m+1 }^--u\le c  |g_{m+1 }|.
\endaligned
$$
For $a^{+}_m+|g_{m}|^{\frac{1}{4}}\le w \le
 a^{-}_{m+1}-|g_{m+1}|^{\frac{3}{5}}$  we have $\ddot w =O(k^{-3})$.
\endproclaim
We start by the last statement. We have $\ddot w =-(\dot w)^3 k''$
with $\dot w\approx 1$ for $a^{+}_m+|g_{m}|^{\frac{1}{4}}\le w \le
 a^{-}_{m+1}-|g_{m+1}|^{\frac{3}{5}}$. In this interval we have
 $ |g_m|^2|u-a_{m}^+|^{-\frac{3}{2}}|u-a_{m}^-|^{-\frac{3}{2}}\le |g_m|^{\frac{5
 }{4}}\ll k^{-3}$ and $ |g_{m+1 }|^2|u-a_{m+1 }^+|^{-\frac{3}{2}}|u-a_{m+1 }^-|^{-\frac{3}{2}}\le |g_{m+1
 }|^{\frac{1
 }{5}}\ll k^{-3}$. We have
$$
 {\ddot E} =2\left ( {\dot u}  \right )^2-2
 u \left (\dot u \right )^3 p''. \tag 1$$
Let us first assume $u\le  a^-_{m +1} -\alpha   $. Then $A(u)
\approx \sqrt{\frac{ |g_m|}{u-a^+_m}} +|g_{m+1}|^2$. By
 (1-3) Lemma 7.3 and by Lemma 7.1 we have
$$
\frac{d^2E}{dp^2}\approx   \frac{1}{(1+ \sqrt{\frac{
|g_m|}{u-a^+_m}} )^2} +    \frac{m}{(1+ \sqrt{\frac{
|g_m|}{u-a^+_m}})^3}\left ( \frac{(a_m^+-a_m^-)^2}{|u-a_m^-|^{\frac
32}|u-a_m^+|^{\frac 32}} + \frac{1}{\langle m \rangle ^3 }\right )
  .
$$
For $|u-a_m^+|\lesssim |g_m|$  we have $ \frac{d^2E}{dp^2}\approx
\frac{m}{|g_m|}.$ For $a^+_{m  }+|g_m| \lesssim u\le   a^-_{m +1}
-\alpha   $ we have
$$\frac{d^2E}{dp^2}\approx 1+       \frac{m|g_m|^2}{ |u-a_m^+|^{3}}   .
$$
Now we consider $u\ge  a^-_{m +1} -\alpha   $. Then $A(u) \approx
\sqrt{\frac{ |g _{m+1}|}{u-a^- _{m+1}}} +|g_{m }|^2$. The two terms
in the right hand side of (1) can be  equal for $|u- a^-
_{m+1}|\approx    |m+1|^{\frac{1}{3}}|g_{m+1 }|^{\frac{2}{3}}$. For
$|u- a^- _{m+1}|\gg   |m+1|^{\frac{1}{3}}|g_{m+1 }|^{\frac{2}{3}}$
we claim that $ \frac{d^2E}{dp^2}\approx 1$. Indeed $A(u)\ll 1$, $
\frac{du}{dp} \approx 1$,   and
$$ \big |     \frac{d^2p}{du^2}      \big | \lesssim \langle u
\rangle ^{-3}+\frac{ m|g_{m+1 }|^2}{|u-a_{m+1 }^-|^{\frac
32}|u-a_{m+1 }^+|^{\frac 32}}\ll 1.$$ For $|u- a^- _{m+1}|\ll
|m+1|^{\frac{1}{3}}|g_{m+1 }|^{\frac{2}{3}}$ we distinguish between
$|u- a^- _{m+1}|\gtrsim |g_{m+1 }|$ and $|u- a^- _{m+1}|\lesssim
|g_{m+1 }|$. If $|g_{m+1 }|\lesssim |u- a^- _{m+1}| \ll
|m+1|^{\frac{1}{3}}|g_{m+1 }|^{\frac{2}{3}}$   then $A(u) \lesssim
1$, $ \frac{du}{dp} \approx 1$ and by (4-5) Lemma 7.3
$$\frac{d^2E}{dp^2}\approx      \frac{-m|g_{m+1 }|^2}{|u-a_{m+1 }^-|^{\frac
32}|u-a_{m+1 }^+|^{\frac 32}}.
$$
For $|g_{m+1 }|\gtrsim |u- a^- _{m+1}|$,  as $u\nearrow a_{m+1 }^-$
then $A(u)$ starts getting larger and $ \frac{du}{dp}$ starts
getting smaller without however matching  $p''$ which is very large,
and we have

$$   u \left (\frac{d u}{dp } \right )^3  \frac{d^2p}{du^2}\approx
(m+1)\frac{|u- a^- _{m+1}|^{\frac{3}{2}  }}{|g_{m+1 }|^{\frac{3}{2}
} } \frac{(a_m^+-a_m^-)^2}{|u-a _{m+1}^-|^{\frac 32}|u-a
_{m+1}^+|^{\frac 32}}\approx \frac{m}{|g_{m+1 }|} $$ and hence
$\frac{d^2E}{dp^2}\approx   -\frac{m+1}{|g_{m+1 }|} $.

\bigskip

\proclaim{Lemma  7.5} Let here $w=u+i0$. There are  fixed $c >0$,
$c_1> 0$ and $n_0$
 such that for $ m \ge n_0$ and for $\alpha \in (1/2,
1)$ then $a_{m+1}^- -  |g _{m+1}|^{\alpha}<u<a_{m+1}^- - c
 |g _{m+1}|   $ implies   $ |\dddot{E}
|\ge c_1 |{m+1}| \, |g_{m+1}|^{4(\frac{1}{2}- \alpha )}.$ Similarly
$a_{m }^+ +  |g _{m }|^{\alpha}>u>a_{m }^++ c
 |g _{m }|   $ implies   $ |\dddot{E}
|\ge c_1 |{m }| \, |g_{m }|^{4(\frac{1}{2}- \alpha )}.$
\endproclaim
We prove the $m+1$ case, the other being similar.  We have
$$\align & \frac{d^3E}{dp^3}=-6\left (\frac{d u}{dp } \right )^4
\frac{d^2p}{du^2}+ 6u\left (\frac{d u}{dp } \right )^5\left (
\frac{d^2p}{du^2} \right )^2- 2
 u \left (\frac{d u}{dp } \right )^4 \frac{d^3p}{du^3}  .\tag 1\endalign$$
 If $c $ is large, in our domain  $\frac{du}{dp}\approx 1$ and
 $    \frac{d^2p}{du^2}
 \approx  \frac{|g _{m+1}|^{2}}{   |u -a^+ _{m+1}|^{\frac{3}{2}}   |u -a^- _{m+1}|^{\frac{3}{2}}          }
   .
 $
 We claim that the dominating term in the rhs of (1) is the third.
We write
$$p^{\prime \prime \prime  }(u)=  \frac 2\pi \sum _n \int _{ g_n}    \frac{
q(t)   }{(t-u)^4}       dt   .$$ For $u$ as in the statement, by (2)
Lemma 7.2
$$\big  |   \frac{d^3p}{du^3}  \big  | \approx
\frac{|g _{m+1}|^{2}}{   |u -a^+ _{m+1}|^{\frac{3}{2}}   |u -a^-
_{m+1}|^{\frac{5}{2}}          }  \gtrsim |g_{m+1}|^{4(\frac{1}{2}-
\alpha )}. $$ Since $\frac{du}{dp}\approx 1$  for $|g _{m+1}|\ll |u
-a^- _{m+1}|$,   we have $ (p^{\prime \prime }) ^2 \ll |p^{\prime
\prime  \prime}| $ and  the last term in the rhs of (1) is the
dominating one. \bigskip

In Lemma 4.6 we have defined $\chi _{int} (k)=\sum _{n\ge n_0} \chi
_1 (\frac{w - a_{n}^+ }{  |g_n|^{ \frac{1}{3}}}  ) \chi _1 (\frac{
a_{n+1}^--w }{ |g _{n+1}|^{ \frac{3}{5}}}  )$ for $n_0\gg 1$. In the
support of $\chi _{int} $ we have $w\approx k$, $  k'= 1+O(\langle
w\rangle ^{-2})$, $k''= O(\langle w\rangle^{-3})$. Now, to complete
the proof of Lemma 4.7 it is enough to prove:

\proclaim{Lemma  7.6} We can extend  the restriction of $p(u)$ on
the support of $\chi _{int} (p(u)) $ to the whole of $\Bbb R$ so
that the extension   $\widetilde{p}(u)$  satisfies the same
relations $w\approx \widetilde{p}$,  $  \widetilde{p}'= 1+O(\langle
u\rangle ^{-2})$, $\widetilde{p}''= O(\langle u\rangle^{-3})$.
\endproclaim
 Recall $p(u)=u+ \frac
1\pi \sum _{n\in \Bbb Z}\int _{  g_n} \frac{ q(t)   }{t-u } dt$ for
$u\in \sigma =\cup \sigma _n$. For $u\in supp \chi _{int} (p(u))$ we
have $p(u)=\widetilde{p}(u)$ with
$$\widetilde{p}(u):=u+ \frac
1\pi \sum _{n<n_0}\int _{  g_n} \frac{ q(t)   }{t-u } dt+ \frac 1\pi
\sum _{n\ge n_0} \chi  _{1} (\frac{u - a_{n}^+ }{  |g_n|^{
\frac{1}{3}}}  ) \chi _{1} (\frac{ a_{n }^--u }{ |g _{n }|^{
\frac{3}{5}}}  ) \int _{ g_n} \frac{ q(t)   }{t-u } dt$$ for a fixed
even smooth cutoff $\chi _{1}$ with $\chi _{1}(t)=0$ for $t$ near 0
and  $\chi _{1}(t)=1$ for $t\gtrsim 1$. If we consider $
\widetilde{p}(u)$ for any $u>a^{+}_{n_0}$ we obtain that
$\widetilde{p}(u)$ satisfies the relations stated in in the
statement. This follows from the fact that near a gap $g _{m}$ we
have
$$\aligned    \widetilde{p}'&= 1+O(\langle
u\rangle ^{-2}) +\frac 1\pi \left ( \chi  _{1} (\frac{u - a_{m}^+ }{
|g_m|^{ \frac{1}{3}}}  ) \chi _{1} (\frac{ a_{m }^--u }{ |g _{m }|^{
\frac{3}{5}}}  ) \int _{ g_m} \frac{ q(t)   }{t-u } dt \right ) '\\&
=1+O(\langle u\rangle ^{-2})+O(|g_m| ^{\frac{4}{5}})\\
 \widetilde{p}'' & = O(\langle u\rangle^{-3})+\frac 1\pi \left ( \chi  _{1} (\frac{u - a_{m}^+ }{  |g_m|^{
\frac{1}{3}}}  ) \chi _{1} (\frac{ a_{m }^--u }{ |g _{m }|^{
\frac{3}{5}}}  ) \int _{ g_m} \frac{ q(t)   }{t-u } dt \right )
''\\& =O(\langle u\rangle^{-3})+  O(|g_m| ^{\frac{1}{5}}).
\endaligned$$

\bigskip

\head \S 8 Estimates on fundamental solutions \endhead

In what follows $\dot f =\frac{d}{dk}f$ and  $  f^\prime
=\frac{d}{dx}f$. Referring to formulas    (2.1) and (2.2),
 we write  $\theta (x,k)$ and $\varphi (x,k)$ for  $\theta (x,w(k))$ and $\varphi
 (x,w(k))$.
 Then we have:

\proclaim{Lemma 8.1} For $x\in [0,1]$ we have
$$\align  &  \big |  \theta (x,k) -\cos (kx)  \big | \le \frac 1{|k|} e^{\frac{x}{|k|}(\|P \| _\infty
+|k^2-E(k)|) }\tag 1
 \\&  \big |  \varphi (x,k) -\frac{\sin (kx)}{k} \big | \le    \frac{1}{|k|^2}\,       e^{\frac{x}{|k|}(\|P \| _\infty
+|k^2-E(k)|) } . \tag 2\endalign
$$
Furthermore there is a fixed constant $C$ such that for $x\in [0,1]$

$$\align  &\big | \dot \theta (x,k) +x \sin (kx)- \int _0^x \theta (s,k)
\frac{\partial}{\partial k}  \frac{\sin \left ( k(x-s) \right )
\left [ P(s)+k^2-E(k) \right ]}{k} ds \big | \le  \frac{C}{\langle k \rangle} \\
& \big |\dot \varphi (x,k) - \frac{\partial}{\partial k} \frac{\sin
(kx)}{k}- \int _0^x \varphi (s,k) \frac{\partial}{\partial k}
\frac{\sin \left ( k(x-s) \right ) \left [ P(s)+k^2-E(k) \right
]}{k} ds \big |  \le \frac{C}{\langle k \rangle ^2} .
\endalign
 $$

\endproclaim

{\it Proof.} The argument is routine.  $\theta (x,k)$ and $ \varphi
(x,k)$ satisfy the following    integral equations:
$$\align & \theta (x,k) =\cos (kx)+\frac 1k \int _0^x \sin \left ( k(x-s)    \right )
\left [      P(s)+k^2-E(k)                       \right ] \theta
(s,k) ds        \tag 3 \\& \varphi (x,k) =\frac{\sin (kx)}{k}+\frac
1k \int _0^x \sin \left ( k(x-s)    \right ) \left [ P(s)+k^2-E(k)
\right ] \varphi (s,k) ds  .  \tag 4  \endalign $$

Now we    write
$$\align  &
\theta (x,k) =\sum _{n=0}^\infty \theta _n(x,k) \, , \quad  \varphi
(x,k) =\sum _{n=0}^\infty \varphi _n(x,k)\\& \theta _0(x,k)=\cos
(kx)\, , \quad  \varphi _0(x,k) =\frac{\sin (kx)}{k} \\& \theta
_{n+1}(x,k) = \frac 1k \int _0^x \sin \left ( k(x-s) \right ) \left
[ P(s)+k^2-E(k)                       \right ] \theta _{n }(s,k)
ds\\& \varphi _{n+1}(x,k) = \frac 1k \int _0^x \sin \left ( k(x-s)
\right ) \left [ P(s)+k^2-E(k) \right ] \varphi _{n }(s,k) ds.
\endalign
$$
 Singling out  $\theta (x,k)$,    we have  for $x_{n+1}=x$: $ \theta _{n+1}(x,k)=$
$$  =\frac 1{k^{n+1}} \int _{0\le x_1 \le \dots \le x_n\le x} \prod
_{j=1}^{n}\left \{
 \sin \left ( k(x_{j+1}-x_{j})    \right ) \left
[ P(x_j)+k^2-E(k)                       \right ] dx_j\right \}
\cos(kx_1).
$$
This implies the following estimate which implies (1):
$$ \big |    \theta _{n+1}(x,k)     \big |
\le   \frac{\left ( \int _0^x(|P(s)|+|k^2-E(k)|) ds \right )
^n}{k^{n+1}\, n!}\le \frac{x^n}{k^{n+1} \, n!} (\|P \| _\infty
+|k^2-E(k)|)^n.
$$
Proceeding similarly we obtain the following inequality, which gives
us (2):
$$ \big |    \varphi _{n+1}(x,k)     \big |
\le   \frac{ \left ( \int _0^x(|P(s)|+|k^2-E(k)|) ds \right )
^n}{k^{n+2}\, n!}\le \frac{x^{n }}{k^{n+2} \, n!} (\|P \| _\infty
+|k^2-E(k)|)^n.
$$

\bigskip
Next we consider
$$\aligned & \dot \theta (x,k) =-x \sin (kx)+ \int _0^x \theta
(s,k) \frac{\partial}{\partial k}  \frac{\sin \left ( k(x-s) \right
) \left [ P(s)+k^2-E(k) \right ]}{k} ds+\\& +\frac 1k \int _0^x \sin
\left ( k(x-s)    \right ) \left [      P(s)+k^2-E(k) \right ] \dot
\theta (s,k) ds
\endaligned $$
and
$$\aligned & \dot \varphi (x,k) = \frac{\partial}{\partial k} \frac{\sin (kx)}{k}+ \int _0^x \varphi
(s,k) \frac{\partial}{\partial k}  \frac{\sin \left ( k(x-s) \right
) \left [ P(s)+k^2-E(k) \right ]}{k} ds+\\& +\frac 1k \int _0^x \sin
\left ( k(x-s)    \right ) \left [      P(s)+k^2-E(k) \right ] \dot
\varphi (s,k) ds.
\endaligned $$
We have for a fixed $C_1$
$$ \aligned &  \big | {\partial _k} \left (  k^{-1}  \sin \left ( k(x-s) \right )  P(s) \right
) +\left ( k^2-E(k)\right ) {\partial _k} \left (  k^{-1}  \sin
\left ( k(x-s) \right ) \right ) + \\& + \left (2k-\dot E (k)\right
)           k^{-1}  \sin \left ( k(x-s) \right ) \big | \le C
\langle k\rangle ^{-1} +C \langle k\rangle ^{-1}\big |2k-\dot E
(k)\big | \le C_1.
\endaligned \tag 8.1
$$
By this estimate, by (1) and  (2) and by the above arguments   we
obtain the last two inequalities of lemma 8.1.

\bigskip

\proclaim{Lemma 8.2} There is a fixed $C$ such that
$$  \align  & \big | \int_0^{1} \theta
^2(x,k) dx-\frac{1}{2}-\frac{\sin (2k)}{4k } \big |\le C k^{-2}
  ; \tag 1 \\&  \big |  k^2 \int_0^{1} \varphi ^2(x,k) dx-
\frac{1}{2}+\frac{\sin (2k)}{4k } \big |\le C k^{-2};\tag 2 \\& \big
|  2k \int_0^{1} \theta  (x,k) \varphi (x,k)dx-\frac{1-\cos
(2k)}{2k} \big |\le C k^{-2}.\tag 3 \endalign  $$
\endproclaim
{\it Proof.} We use the notation in the proof of Lemma 8.1. To prove
(1) is enough to show that $ \int_0^{1} \cos (x k)\theta _1 (x,k)
dx=O(k^{-2}).$ By its definition  and elementary computation
$$\theta _1 (x,k)=\frac{\sin (x k)}{2k} \left [ \int _0^xP(s) ds+(k^2-E
)x\right ] +O(k^{-2}).\tag 4$$ By $k^2-E=-2Q_0+o(1)$ and elementary
integration, (1) follows.

 To prove
(2) is enough to show that $ \int_0^{1} \sin (x k)k\varphi _1 (x,k)
dx=O(k^{-2}).$ By its definition and elementary computation
$$k\, \varphi _1 (x,k)=-\frac{\cos (x k)}{2k} \left [ \int _0^xP(s) ds+(k^2-E
)x\right ] +O(k^{-2}) .\tag 5$$ Elementary integration gives (2).

To prove (3) is enough   $ \int_0^{1} \left ( \cos (x k)k\varphi _1
(x,k) + \sin (x k)\theta _1 (x,k) \right ) dx=O(k^{-2}).$   So (3)
follows from $ \int_0^{1} \left ( \cos (x k)k\varphi _1 (x,k) + \sin
(x k)\theta _1 (x,k) \right ) dx=$ $$= - \frac{1}{2k} \int_0^{1}
\cos (2kx) \left [ \int _0^xP(s) ds+(k^2-E )x\right ] dx+O(k^{-2}).
$$

\bigskip

\head \S 9 Proof  of Lemma 4.4 \endhead

Formulas (2.2), (2.3), (2.7) and Lemma 8.1 imply with Lemma 9.1
below gives the first claim of Lemma 4.4.

\proclaim{Lemma 9.1} There is a fixed constant $C>0$ such that for
$w=u+iv$ with $v=0$ and for  $a_n^+ +Cn^5|g_n|\le u\le  a_{n+1}^{-}
-Cn^5|g_{n+1}|$, we have
$$ \big | N^2(k)-1 \big |+ \big | \frac{ \sin k }{k \varphi (k) } -1\big |+  \big |\frac{\varphi ^\prime (k)  - \theta (k) }{2k\varphi (k)}
\big |\le  \frac{C}{\langle k \rangle }. $$
 \endproclaim

{\it Proof}. Expanding in the definition of $N^2(k)$ we write
$$\aligned &N^2(k)= \int_0^{1} \left ( \theta (x,k) +  \frac{\varphi ^\prime (k)  - \theta (k)
}{2\varphi (k)} \varphi   (x,k) \right ) ^2  dx  +  \int_0^{1}
   \frac{ \sin ^2k }{ \varphi ^2(k) }  \varphi  ^2 (x,k)   dx=\\&  = \int_0^{1}
   \left [   \theta ^2(x,k)  +2\theta (x,k)  \varphi   (x,k)    \frac{\varphi ^\prime (k)  - \theta (k)
}{2\varphi (k)}+     \left (\frac{\varphi ^\prime (k)  - \theta (k)
}{2\varphi (k)} \varphi   (x,k)\right )^2        \right ] +\\&
+\int_0^{1}
   \frac{ \sin ^2k }{ \varphi ^2(k) }  \varphi  ^2 (x,k)   dx
.\endaligned$$  We recall now from formulas (1.4) and (3.1) in
\cite{F2}: $$\dot E(k) =\frac{2\sin k }{ \varphi (k) N^2(k)}. \tag
9.1$$  For $a_n^+ +Cn^5|g_n|\le u \le a_{n+1}^{-} -Cn^5|g_{n+1}|$
we have $\big | \dot E(k)-2k\big | \lesssim \langle n \rangle ^{-2}
$ by claim 4 in the statement of Lemma 7.1. So $ N^2(k) =\frac{\sin
k }{k \varphi (k)  }    \left ( 1+O_2(\frac 1{k^2}) \right )     $
with $O_2$  a big $O$.  Hence by Lemma 8.2 we obtain, for a certain
number of big O's,

$$  \aligned & \left ( \frac{ \sin k }{k \varphi (k) } \right )
^2-2\frac{\sin k }{k \varphi (k)  }      \frac{1+O_2(k^{-2})}{1
-\frac{\sin (2k)}{2k}          +O _1( k^{-2}) } +\left (
\frac{\varphi ^\prime (k) - \theta (k) }{2k\varphi (k)} \right ) ^2
+\\& +\frac{1 +\frac{\sin (2k)}{2k}          +O _3( k^{-2}) }{1
-\frac{\sin (2k)}{2k} +O_1 ( k^{-2}) } +O _4( k^{-1}) \frac{\varphi
^\prime (k) - \theta (k) }{2k\varphi (k)} =0\endaligned
$$
which implies $  \frac{ \sin k }{ k\varphi (k) }  =
\frac{1+O_2(k^{-2})}{1 -\frac{\sin (2k)}{2k}          +O _1( k^{-2})
}\pm \sqrt{\Delta}$,  for $$
  \Delta= {
\frac{\frac{\sin ^2(2k)}{4k^2}}{ \left ( 1- \frac{\sin (2k)}{2k}
\right )^2}+ O(k^{-2})- \left ( \frac{\varphi ^\prime (k)  - \theta
(k) }{2k\varphi (k)} \right ) ^2  -O _4( k^{-1}) \frac{\varphi
^\prime (k) - \theta (k) }{2k\varphi (k)}}.
$$
So $ \frac{ \sin k }{k \varphi (k) } =1+O(\frac 1k)$ and $
\frac{\varphi ^\prime (k)  - \theta (k) }{2k\varphi (k)} =O(\frac 1{
{k}}) $ when $k\in \Bbb R$ and under the restriction $a_n^+
+Cn^5|g_n|\le u\le  a_{n+1}^{-} -Cn^5|g_{n+1}|$.

\bigskip

\noindent The proof of the second claim of Lemma 4.4 is trickier and
proceeds in several steps. We  set $\|m_\pm ^0( k)\| _2=(\int
_0^1|m_\pm ^0(x,k)| ^2dx) ^{\frac{1}{2}}$. First of all we consider
the Fourier series expansion $m_{\pm }^{0}(x,k)=\sum
\widehat{m}_{\pm}(\ell ,k)e^{2\pi i \ell x}$ and show that the $L^2$
is concentrated in  two harmonics. One harmonic is $|\widehat{m_\pm
}(0,k )|\approx \| {m_\pm }(k )\| _2 $. If there is an $n$ such that
$|n\pi +k|\ll 1$ then also $ \widehat{m_\pm }(\mp n,k )$ can be
significant. We then bound $N^ {-1}(k)$, $m^{\pm }(k)k^{-1}N^
{-1}(k)$ in terms $\| {m_\pm }(k )\| _2 $. Next, we express $m_\pm
^0(x,k)$ in terms of $\varphi (x,k)$ and $\theta (x,k)$, we expand
the latter in terms of $\sin(kx)$ and $\cos(kx)$ and a reminder, and
we conclude that in $L^\infty$ sense $m_\pm ^0(x,k)$ can be
approximated by the two terms of the Fourier expansion discussed
above. Next we look at the normalization of the Bloch functions.
Slightly off the slits we have $|\widehat{m_\pm }(0,k )|\gg
|\widehat{m_\pm }(\mp n,k )|$ and so $1\approx \widehat{m_+ }(0,k )
\widehat{m_- }(0,k )$. From this we get the desired bound on $|{m_+
}(x,k )  {m_- }(y,k )|\lesssim 1$ off the slits. Near the slits we
have $|\widehat{m_\pm }(0,k )|\approx |\widehat{m_\pm }(\mp n,k )|$
so to exploit the normalization we have to exclude a significant
cancelation in a certain formula. Let us start with the first step,
and show that there are at most two significant harmonics.

\proclaim{Lemma 9.2} For all $n$ except  possibly for  one $ n_0$,
we have $|n\pi +k|\gtrsim 1$. Then for $n\neq   0, n_0$   we have $
\sum _{n \neq 0,n_0 }
  | \widehat{m_+}(n  ) |^2 \le C
|k| ^{-2}\| m_+ ^0( k)\| _2^2$    for a fixed $C$.  Furthermore for
a fixed $C$ we have $ \sum _{n \neq 0,n_0 } n^2
  | \widehat{m_+}(n  ) |^2 \le C
 \| m_+ ^0( k)\| _2^2$. The same statement holds for $m_-
^0$ with $ n_0$ replaced  by $-n_0$. If for all $n\neq 0$ we have
$|n\pi +k|\gtrsim 1$ we can extend the above inequalities to the sum
on  all $n \neq 0$.
\endproclaim
{\it Proof.} The first sentence is straightforward. We will assume
there is  $n_0$ with $|n_0\pi +k|\ll 1$. If such $n_0$ does not
exist, the proof is almost the same. Set for $n\neq 0, n_0$

$$ \aligned & \widehat{m_+}(n)+\sum _{\ell \neq  0, n_0  } T(n,\ell
)  \widehat{m_+}(\ell ) =-\frac{ \widehat{P}( n) \widehat{m_+}(0 ) +
\widehat{P}( n-n_0) \widehat{m_+}(n_0 ) }{ 4\pi n (n\pi +k) } \\&
T(n,\ell ) =\frac{\widehat{P}(0) +k^2-E }{4\pi n (n\pi +k)} \delta
(n-\ell )+\frac {\widehat{P}(n-\ell )  }{ 4\pi n (n\pi +k) }
.\endaligned \tag 1
$$
We have $\widehat{P}(0) +k^2-E=O(k^{-2}) $ by Lemma 5.2 and by
$\widehat{P}(0)=2Q_0$.  Equation (1) is of the form $(I+T)u=f$ where
$\| f\| _{l^2} \le C |k|^{-1} \| m_\pm ^0( k)\| _2 $ and $\| f\|
_{l^2_1} \le C   \| m_\pm ^0( k)\| _2 $ where $\| f\| _{l^2_1}^2=
\sum _{n \neq  0, n_0  } n^2|f (n)|^2.$ By
$$ \sup _\ell \sum _n \left ( |k T(n,\ell
) | +|n T(n,\ell ) |\right ) +\sup _n \sum _\ell  \left ( |k
T(n,\ell ) | +|n T(n,\ell ) |\right ) \le C
 $$  for $u\in l^2 (\Bbb N)$ we have $\sum _n |(Tu)(n)|^2 \le
C |k|^{-2} \| u\| _{l^2}^{2}$ and $\sum _n n^2 |(Tu)(n)|^2 \le C \|
u\| _{l^2}^{2}$. So inverting and after a Neumann expansion, we see
that for $n\neq 0, n_0$ we have
$$\widehat{m_+}(n)=-\frac{ \widehat{P}( n) \widehat{m_+}(0 ) +
\widehat{P}( n-n_0) \widehat{m_+}(n_0 ) }{ 4\pi n (n\pi +k) } +
\widehat{e} (n)\tag 2$$ where $\widehat{e}=\sum _{m=1}^\infty
(-)^mT^mf$ satisfies  $ \widehat{e} (0)=\widehat{e}(n_0)=0$ and
$$\| \widehat{e}  \| _{l^2 } \le \| \widehat{e}  \| _{l^2_1} \le C
\| \sum _{m=0}^\infty (-)^mT^mf\| _{l^2} \le C |k|^{-1} \sum
_{m=0}^\infty (C|k|^{-1})^m   \| m_+ ^0( k)\| _2  .$$ Hence by (2),
$ \sum _{n \neq 0, n_0 }
 |\widehat{m_+} (n)|^2\le C |k|^{-2}\| m_+
^0( k)\| _2^2$ and  $ \sum _{n \neq 0, n_0 }n^2
 |\widehat{m_+} (n)|^2\le C  \| m_+
^0( k)\| _2^2$. The proof for $m_-$ is similar.

\bigskip

We express now the Bloch functions in terms of the fundamental
solutions as in \S 2. Using   the notation in the proof of Lemma 8.1
and for $m^\pm (k)= \frac{\varphi ^\prime (k)  - \theta (k)
}{2\varphi (k)} \pm i \frac{\sin k}{\varphi  (k)}  $,
$$\aligned   m_\pm  ^0(x ,k)    &=e^{\mp ikx}  \left ( \frac{\cos(x k)}{N} + \frac{m^\pm (k)}{kN} \sin(x k)\right )
  +\\& + \frac{1}{N }e^{\mp ikx}\sum _{j=1}^\infty
 \theta _j(x,k)  + \frac{m^\pm (k)}{kN}  e^{\mp ikx}\sum
_{j=1}^\infty k\varphi _j (x,k)  .
\endaligned \tag 9.2$$
  By the proof of Lemma 8.1  $ \big | e^{\mp ikx}\sum _{j=1}^\infty
 \theta _j(x,k)\big | +\big |  e^{\mp ikx}\sum _{j=1}^\infty
 k\, \varphi _j(x,k)\big | =O(k^{-1})$ for $x\in [0,1]$.  Then we have

$$\aligned &m_\pm  ^0(x ,k)    =A_\pm (k)    + B_\pm (k)    e^{\mp i
2kx} + O(N^{-1}k^{-1}) + O(m^\pm (k)N^{-1}k^{-2}),
\\& A_\pm (k) =\frac{1 -i\, k^{-1}  m^\pm (k)  }{2N(k)}=\frac{1+ \frac{\sin ( k)}{ k\varphi (k) }  \mp i
\frac{\varphi ^\prime (k)  - \theta (k) }{2k\varphi (k)} }{2N(k)}\\&
   B_\pm (k)= \frac{1+i\, k^{-1}  m^\pm (k) }{2N(k)}= \frac{1- \frac{\sin ( k)}{ k\varphi (k) }  \pm
i \frac{\varphi ^\prime (k)  - \theta (k) }{2k\varphi (k)}  }{2N(k)}
.\endaligned
$$
We have
$$\widehat{m_\pm  ^0(\cdot  ,k)}(0)=A_\pm (k)+ O\left (k^{-1} B_\pm (k) \right )+ O(N^{-1}k^{-1})
+ O(m^\pm (k)N^{-1}k^{-2}).\tag 9.3 $$ For the $n_0$ of Lemma 9.2 we
have
$$\widehat{m_\pm  ^0(\cdot  ,k)}(\pm n_0)= B_\pm (k)  \frac{ e^{\pm 2i ( \pi n_0-k)} -1}{\pm 2i ( \pi n_0-k)}
+ O(N^{-1}k^{-1}) + O(m^\pm (k)N^{-1}k^{-2}).\tag 9.4
$$

\proclaim{Lemma 9.3} For  a  fixed $C>0$
$$\align  & |1/N(k)| +  |\frac{   m^\pm (k) }{ kN(k)} |  +|A_{\pm }| + |B_{\pm }|\le C \| m_\pm ^0(
k)\| _2
 . \endalign$$
\endproclaim
\noindent {\it Proof.} Take $$A_\pm (k)=\frac{1}{2N}+i\frac{ m^\pm
(k)}{2kN}\, , \quad B_\pm (k)=\frac{1}{2N}-i\frac{ m^\pm (k)}{2kN}.
\tag 1 $$ Then by (9,3-4)
$$ \align & |A_\pm (k)| \lesssim   \| m_\pm ^0( k)\| _2+ O(
N^{-1}k^{-1})+O(m^\pm  (k)N^{-1}k^{-2})  \tag 2\\& |B_\pm (k)|
\lesssim  \| m_\pm ^0( k)\| _2+  O( N^{-1}k^{-1})+O(m^\pm
(k)N^{-1}k^{-2}) .\tag 3 \endalign$$   By the triangular inequality,

$$|1/N(k)| \le |A_\pm (k)|+|B_\pm (k)| \lesssim \| m_\pm ^0( k)\| _2+ O(m^\pm (k)N^{-1}k^{-2}).\tag 4$$
  By (1-3)  if one of $|1/N(k)|$ and $|m^\pm
(k)N^{-1}k^{-1}| $ is $\gg \| m_\pm ^0( k)\| _2$, then
$|1/N(k)|\approx |m^+ (k)N^{-1}k^{-1}| $ . Then $|1/N(k)|\lesssim \|
m_\pm ^0( k)\| _2$ by (4) and   $|m^\pm (k)N^{-1}k^{-1}|\lesssim \|
m_\pm ^0( k)\| _2 .$ So   these last two formulas hold. So $|A_\pm
(k)|+|B_\pm (k)| \lesssim \| m_\pm ^0( k)\| _2 $ by (2-3) .

\bigskip

We now  use the normalization of Bloch functions $1=\int _0^1m_+
^0(x, k)m_- ^0(x, k)dx.$ We will denote $n_0$ by $n$. By (9,3-4) and
Lemma 9.3
$$1=\left [ A_+ (k)A_- (k)+B_+ (k)B_- (k)\right ]+ O(k^{-1}\| m_+
^0( k)\| _2\| m_- ^0( k)\| _2).\tag 9.5$$ We can also write by Lemma
9.2
$$1= \widehat{m} _{+} (0,k ) \widehat{m} _{-} (0,k )+
\widehat{m} _{+} (- n ,k )\widehat{m} _{-} (  n ,k )+ O(k^{-2}\| m_+
^0( k)\| _2\| m_- ^0( k)\| _2).\tag 9.6$$

\proclaim{Lemma 9.4} Suppose that $k$ is in a region such that
$|\widehat{m} _{\pm } (\mp n ,k )|\le \frac{1}{2}|\widehat{m} _{\pm
} (0,k )|$. Then for a fixed constant $C$ we have $ \| m_+ ^0( k)\|
_2 \| m_-^0( k)\| _2\le C $. As a consequence also $| m_+ ^0(x, k)
m_- ^0(y, k)|<C$ for some fixed $C$.
\endproclaim  {\it Proof.} $ \| m_+ ^0( k)\|
_2 \| m_-^0( k)\| _2\le C_1 $ for a fixed $C_1$ by $|\widehat{m}
_{\pm } (0,k )|\approx \| m_\pm  ^0( k)\| _2$ and (9.6). By (9.2-4)
and Lemma 9.3 we get $| m_+ ^0(x, k) m_- ^0(y, k)|<C$ for a fixed
$C$.

\bigskip
Lemma 9.4 applies to the case when $k$ is not close to the slits,
for example if there is no $n$ with $| \pi n+k|\ll 1$. Let us
suppose $n=n_0$ exists. We have:

\proclaim{Lemma 9.5} Consider $k=p+iq$ and corresponding $w=u+iv$.
For a fixed constant $C$  and for $I_n(u,v )=\frac{v}{\pi}\int
_{a^{-} _{n}}^{a^{+} _{n}} \frac{q(t)\, dt}{(t-u)^2+v^2}   $ we have
$$ \big |  |\widehat{m} _{\pm} (\mp n,k )|^2-\frac 12 \,  \frac{|u|}{ |n |\pi} \frac{|I_{n }(u,v)|}{|q|}
\| m_\pm ^0( k)\| _2^2\big | \le C |k |^{-2}  \| m_\pm ^0( k)\|
_2^2.\tag 9.7$$
\endproclaim
{\it Proof.} It is enough to consider $q>0$. Set $m=m_\pm ^{0}$.
Multiply by $\overline{m }$ the equation $  m'' \pm 2ikm'
-P(x)m+(E-k^2)m=0  $, integrate in $[0,1]$ and take imaginary part
to obtain
$$
2q \, \Im \int_0^1 \overline{m} \, m' dx = \Im (E-k^2) \| m  \|
_2^2.
   \tag 9.8 $$
Then by Lemma 9.2, $\int_0^1 \overline{m} \, m' dx= -2\pi n i
|\widehat{m} (\mp n,k)|^2+O(k^{-1}\| m\| _2^2)$ and $\Im (E-k^2)
=(v-q)(u+p)+(u-p)(v+q) =-(2u +O(u^{-1})) I_n+O((|v|+|q|) u^{-1}).$
Substitute in (9.8), divide by $4\pi q|n|$ and use $0 \le v \le q$
and $ I_n\le q$.

\bigskip

\proclaim{Lemma 9.6} There is a fixed $\Gamma >0$ such that for
$|p-n\pi |>\Gamma |g_n|$ we have $|\widehat{m} _{\pm } (\mp n ,k
)|\le \frac{1}{2}|\widehat{m} _{\pm } (0,k )|$.
\endproclaim
{\it Proof.} For $\Gamma \gg 1$ we have $I_n/q\ll 1$ in (9.7) by
Lemma 6.2. Furthermore $ {|u|}/{ (|n |\pi )}= 1+O(1/|u|).$

\bigskip

\proclaim{Lemma 9.7} There is a fixed $C$ such that for any  fixed
$\Gamma
>0$ there   are $N,\delta
>0$ such that for $|n|\ge N$  in the region $|p-n\pi |<\Gamma |g_n|$, and $0\le q \le
\delta |g_n| $ we have $\| m_+ ^0( k)\| _2\| m_- ^0( k)\| _2\le C.$
\endproclaim
{\it Proof}. First of all, by Lemma 9.5 for $k$ near the slit we
have $| \widehat{m}_\pm  ^0(j,  k) |^2 \approx \| m_+ ^0( k)\|
_2^2/2$ for $j=0,\mp n$. These harmonics could be  large with a
  large cancelation in (9.6). We have
$$\aligned &    A_+ (k)A_- (k)+B_+ (k)B_- (k) =\frac{\left (1+ \frac{\sin ( k)}{
k\varphi (k) } \right )^2   +\left ( \frac{\varphi ^\prime (k)  -
\theta (k) }{2k\varphi (k)} \right )^2 }{4N^2(k)} +\\& +\frac{\left
(1- \frac{\sin ( k)}{ k\varphi (k) } \right )^2 +\left (
\frac{\varphi ^\prime (k)  - \theta (k) }{2k\varphi (k)} \right )^2
}{4N^2(k)}  =\frac{1+ \left (  \frac{\sin ( k)}{ k\varphi (k) }
\right )^2 +\left ( \frac{\varphi ^\prime (k)  - \theta (k)
}{2k\varphi (k)} \right )^2 }{2N^2(k)} . \endaligned $$ By  (9.1) we
have
$$A_+ (k)A_- (k)-B_+ (k)B_- (k)=\frac{\sin ( k)}{
k\varphi (k) N^2}=
 \frac{\dot E}{2k}= \frac{ w}{ k \frac{dk}{dw}}\approx
\frac{ 1}{   \frac{dk}{dw}}.\tag 9.9$$ We have $ \| m_+ ^0( k)\| _2
\| m_- ^0( k)\| _2\approx $ $$ \approx |A_+ (k)A_- (k)|\le
1+2/|k'(w)|+O(k^{-1}\| m_+ ^0( k)\| _2 \| m_- ^0( k)\| _2).$$
 By  Lemma 9.8  below we obtain  $ \| m_+ ^0( k)\| _2 \| m_- ^0( k)\| _2\lesssim 1.$

\proclaim{Lemma 9.8} There is a fixed $C>0$ such that for any  fixed
$\Gamma
>0$ there are $N,\delta
>0$ such that for $|n|\ge N$     in the region $|p-n\pi |<\Gamma |g_n|$, and $0\le q \le
\delta |g_n| $ we have $|k'(w) |  > C.$
\endproclaim
{\it Proof.} We recall by Schwartz Christoffel formula,   for $c_m
\in g_m $ with $k'(c_m )=0$,

$$k'(w)= a\frac{1-w/c_n }{\sqrt{ (1-w/a^{-}_n ) (1-w/a^{+}_n)    }     }
\prod _{\ell \neq n} \frac{1-w/c_\ell }{\sqrt{ (1-w/a^{-}_\ell )
(1-w/a^{+}_\ell)    }     } , \tag 1$$ see (1.10) \cite{MO}. By
taking the derivative $dq/du$ in (5.1) we see that $q'(u)>0$ for
$u\in ( a_m ^-, a_m ^-+ \varepsilon |g_m |] $ and $q'(u)<0$ for
$u\in [ a_m ^+- \varepsilon |g_m |, a_m ^+) $ for a fixed
$\varepsilon
>0$. So
$c_m \in [   a_m ^-+ \varepsilon |g_m|,a_m ^+- \varepsilon |g_m | ]$
for a fixed sufficiently small $\varepsilon
>0$. Since the infinite product in (1) has value approximately 1 for $w$ near $\pi n$,
because of the second factor   there is a fixed $C>0$ such that
$|k'(w) |  > C $ holds for $ w=u+iv$ with $|v|\le \delta |g_n|$ and
either  $a^{-}_n -\Gamma |g_n|\le u\le  a_n^-+ \varepsilon |g_n|/2$
or $a_n^+- \varepsilon |g_n|/2\le u \le a^{+}_n +\Gamma |g_n| $ with
$\delta
>0$ fixed and sufficiently small. Now we need to show that the
values  of $k=p+iq$ in the statement are inside this region in the
$w$ plane. First of all $0\le q \le \delta |g_n| $ implies $0\le
v\le q \le \delta |g_n| $ by (6.1). For $0\le v  \le \delta |g_n| $
and $a_n^-+ \delta |g_n|\ll u \ll a_n^+- \delta |g_n|$ by Lemma 6.3
we have $q\approx I_n(u,v)$ and by Lemma 6.4 we have
$I_n(u,v)\approx |g_n|^{ \frac{1}{2}} \sqrt{\min   \{   |u-  a^{-}
_{n}|,    |u- a^{+} _{n}| \}} \gg \delta |g_n|$. Obviously the
latter is incompatible with $q\le \delta |g_n|$. So for small
$\delta >0$ we have either $a^{-}_n -\Gamma |g_n|\lesssim u\le
a_n^-+ \varepsilon |g_n|/2$ or $a_n^+- \varepsilon |g_n|/2\le u
\lesssim a^{+}_n +\Gamma |g_n| $.

\bigskip Notice that by (9.9) and (1) in Lemma 9.8 we can see that
near the tip of the slit  the product $\| m_+^0(k)\| _2\| m_-^0(k)\|
_2$ is unbounded.

\bigskip \head \S 10 Proof of Lemma 4.5: case    $u\in [a^{+}_{n}+
|g_n|^{\frac{1}{4}}, a^{-}_{n+1}-  |g
_{n+1}|^{\frac{3}{5}}]$\endhead

Lemma 4.5 consists in 3 claims. The third one follows immediately
from Lemma 4.4 by the Cauchy integral formula. The first two claims
follow immediately from the Cauchy integral formula from Lemma 10.1
which is an improvement of the second claim of Lemma 4.4 in the case
when $a^{+}_{n}+ |g_n|^{\frac{1}{4}}\le u \le a^{-}_{n+1}-  |g
_{n+1}|^{\frac{3}{5}}$. We will assume this restriction on $u$
everywhere below in this section:

\proclaim{Lemma 10.1} For $a^{+}_{n}+ |g_n|^{\frac{1}{4}}\le u \le
a^{-}_{n+1}-  |g _{n+1}|^{\frac{3}{5}}$ with $n\ge n_0$ and for
$1\ge v\ge 0$, there is a $C$ such that $ \big | m_- ^0(x ,k)m_+
^0(y ,k)-1 \big | \le \frac{C}{|k|} .$
\endproclaim

Set $m_\pm (x)= m_\pm  ^0(x ,k)$ and   consider the expansion $m_\pm
(x)= \sum _n e^{2\pi i n x}\widehat{m_\pm }(n).$ For $a^{+}_{n}+
|g_n|^{\frac{1}{4}}\le u \le a^{-}_{n+1}-  |g _{n+1}|^{\frac{3}{5}}$
we have the following strengthening  of Lemma 9.3:

\proclaim{Lemma 10.2 } We have  for a fixed $C$ $$||\widehat{m_\pm
}(0 ) |^2 -\| m_\pm ^0( k)\| _2^2|  +\sum _{n \neq 0 }
  | \widehat{m_\pm}(n  ) |^2 \le C
|k| ^{-2}\| m_\pm ^0( k)\| _2^2.$$
\endproclaim
\noindent {\it Proof .}   By (9.8) and   Lemma 6.2 we get
 $  \Im \int _0^1 \bar m _\pm m^\prime _\pm dx = O( u^{-2} )\| m_\pm
^0( k)\| _2 ^2.$ By Lemma 9.2 $\Im \int _0^1 \bar m _\pm m^\prime
_\pm dx =\mp 2\pi i n_{0}|\widehat{m}_\pm (\mp
n_{0},k)|^2+O(k^{-1}\| m_\pm ^0( k)\| _2 ^2)$. Hence
 $| \widehat{m_\pm}(\pm n_0  ) |  \le C |k| ^{-1}\| m_\pm
^0( k)\| _2  $. The latter and Lemma 9.2 imply the inequality
$||\widehat{m_\pm }(0 ) |^2 -\| m_\pm ^0( k)\| _2^2  |\le C |k|
^{-2}\| m_\pm ^0( k)\| _2^2.$

\bigskip

Now we have the following   lemma:

\proclaim{Lemma 10.3} For   fixed constants we have $| N
^{-1}|\approx  |m^\pm (k)N^{-1}k^{-1}|\approx |A_{\pm }| \approx \|
m_\pm ^0( k)\| _2 $ and   $|B_{\pm }| \lesssim |k|^{-1} \| m_\pm ^0(
k)\| _2$.
\endproclaim
\noindent {\it Proof of Lemma 10.3.} We get $|B_{\pm }| \lesssim
|k|^{-1} \| m_\pm ^0( k)\| _2$ by Lemmas 10.2 \& 9.3 and by (9.4) .
By $|\widehat{m}_\pm (0,k)|\approx \| m_\pm ^0( k)\| _2 $ and by
(9.3) we get $|A_{\pm }| \approx \| m_\pm ^0( k)\| _2 $. By
definition of $A_{\pm }$ and of $B_{\pm }$, estimates $|B_{\pm }|
\lesssim |k|^{-1} \| m_\pm ^0( k)\| _2$ and $|A_{\pm }|\approx \|
m_\pm ^0( k)\| _2 $  imply $|1/N(k)|\approx |m^\pm
(k)N^{-1}k^{-1}|\approx \| m_\pm ^0( k)\| _2$.

\bigskip

\proclaim{Lemma 10.4} We have for   fixed constants $\| m_\pm ^0(
k)\| _2  \approx 1$.\endproclaim

\noindent {\it Proof of Lemma 10.4.} By Lemma 10.3   $|A_\pm
(k)|\approx |1/N|\approx \| m_\pm ^0( k)\| _2 $. By  Lemma 10.2 and
the formulas immediately above (9.5)
  $ m_\pm ^0( x,k) \approx A_\pm (k)$. In particular
by  normalization of Bloch functions we have the following which
gives us Lemma 10.4:
$$1=\int _0^1m_+ ^0( x,k)m_- ^0( x,k)dx\approx A_+ (k) A_- (k).$$

\bigskip
{\it Proof of Lemma 10.1. } We write

 $$ m_- ^0(x ,k)m_+ ^0(y ,k)=\frac{\left (1+ \frac{\sin ( k)}{
k\varphi (k) } \right )^2 +\left ( \frac{\varphi ^\prime (k)  -
\theta (k) }{2k\varphi (k)} \right )^2 }{4N^2(k)} + O(k^{-1}) .$$ By
the expansion of $N^2(k)$ in Lemma 9.1 and by Lemmas 8.2,   10.3 and
10.4, we have

$$2N^2(k)=1+\left ( \frac{\sin ( k)}{
k\varphi (k) } \right )^2 +\left ( \frac{\varphi ^\prime (k)  -
\theta (k) }{2k\varphi (k)} \right )^2+O(k^{-1}).$$ This, formula
(9.1) and $\frac{\dot E}{2k} = 1+ O(k^{-1})$,   Lemma 7.1 (5), imply
the following,
$$\frac{1+ \left (  \frac{\sin ( k)}{
k\varphi (k) } \right )^2 +\left ( \frac{\varphi ^\prime (k)  -
\theta (k) }{2k\varphi (k)} \right )^2 }{4N^2(k)}   +
\frac{\frac{\sin ( k)}{ k\varphi (k) }}{2N^2(k)}
=\frac{1}{2}+\frac{\dot E}{4k}+ O(k^{-1})= 1+ O(k^{-1}),  $$
 which ends the proof of Lemma 10.1.  An immediate consequence of
Lemma 10.1, of the geometry of the comb $\Cal K$ and of the Cauchy
integral formula is:

\proclaim{Lemma 10.5} There is $n_0$ such that if for $n\ge n_0$ we
have $a^{+}_{n}+  |g_n|^{\frac{1}{4}}\le u \le
\frac{a^{+}_{n}+a^{-}_{n+1}}{2} $ and $v=0$,  then there is a $C$
such that for the corresponding $k=p+i0$ we have $  | \partial _k
(m_- ^0(x ,k)  m_+ ^0(y ,k))\big | \le \frac{C}{k |k-\pi n|}
 .$ If  $
 \frac{a^{+}_{n}+a^{-}_{n+1}}{2} \le u \le a^{-}_{n+1}-
 |g _{n+1}|^{\frac{3}{5}}$ then $  | \partial _k
(m_- ^0(x ,k)  m_+ ^0(y ,k)) \big | \le \frac{C}{k |k-\pi (n+1)|}
 .$
\endproclaim

\enddocument